\documentclass[12pt]{amsart}
\usepackage{graphicx} 
\usepackage{mathtools}
\usepackage{amsmath}
\usepackage{amssymb}
\usepackage{amsthm}
\usepackage{circuitikz}
\usepackage{dirtytalk}
\newcommand{\diam}{\text{diam}}

\usepackage[maxbibnames=99,maxcitenames=99,
backend=biber,
style=alphabetic,
sorting=nyt
]{biblatex}
\usepackage{amsmath}
\usepackage{amssymb}
\usepackage{amsthm}
\usepackage{dsfont}
\usepackage{mathrsfs}
\usepackage[all]{xy}
\usepackage[mathcal]{eucal}
\usepackage{verbatim}  
\usepackage{fullpage}  
\usepackage{hyperref}
\usepackage[shortlabels]{enumitem}
\usepackage{upgreek}
\usepackage{bm}
\usepackage{color,soul}
\usepackage{graphicx}
\usepackage{commath}
\usepackage{float}

\newcommand{\cC}{\mathcal{C}}
\newcommand{\cK}{\mathcal{K}}

\newcommand{\Hs}{\mathcal{H}^s}
\newcommand{\CXY}{\mathcal{C}(X,Y)}

\newcommand{\su}{\subseteq}

\newcommand{\smr}{\preccurlyeq}

\newcommand{\RR}{\mathbb{R}}
\newcommand{\R}{\mathbb{R}}

\newcommand{\NN}{\mathbb{N}}

\renewcommand*{\d}{\mathop{\kern0pt\text{\upshape d}}\!{}}

\renewcommand{\epsilon}{\varepsilon}

\newtheorem{theorem}{Theorem}[section]
\newtheorem{proposition}[theorem]{Proposition}
\newtheorem{lemma}[theorem]{Lemma}
\newtheorem{corollary}[theorem]{Corollary}

\theoremstyle{definition}
\newtheorem{definition}[theorem]{Definition}
\newtheorem{remark}[theorem]{Remark}

\newtheorem{notation}{Notation}

\usepackage{chngcntr}
\counterwithout{equation}{section}

\renewcommand*{\d}{\mathop{\kern0pt\text{\upshape d}}\!{}}

\newcommand{\nor}[1]{\left\lVert #1 \right\rVert}

\definecolor{red}{rgb}{1.0, 0.33, 0.64}

\PassOptionsToPackage{hyphens}{url}\usepackage{hyperref}
\usepackage{xurl}
\usepackage{esint}

\newtheorem{conjecture}{Conjecture}

\hypersetup{
    colorlinks=true,
    linkcolor=blue,
    filecolor=magenta,      
    urlcolor=cyan,
    pdftitle={Typical Kakeya-Type Sets and Their Typical Code Sets},
    }

\title{Typical Kakeya-type Sets and Their Typical Code Sets}

\author{Josef Bell}
\address{Department of Mathematics, Yale University, New Haven, CT 06511, USA}
\email{josef.bell@yale.edu}

\author{Ezekiel Fleissner-Kates}
\address{Department of Mathematics, University of Chicago, Chicago, IL 60637, USA}
\email{eofleissnerkates@uchicago.edu}

\author{Tamás Keleti}
\address
{Institute of Mathematics, ELTE E\"otv\"os Lor\'and University, 
P\'az\-m\'any P\'e\-ter s\'et\'any 1/c, H-1117 Budapest, Hungary and
HUN-REN Alfr\'ed R\'enyi Institute of Mathematics,
Re\'altanoda utca 13-15, H-1053 Budapest Hungary}
\email{tamas.keleti@gmail.com}

\thanks{The authors would like to thank the Budapest Semesters in Mathematics program
for providing the framework under which this research was conducted.
The third author was also supported  
by the Hungarian National Research, Development and Innovation Office -- NKFIH, 146922 and 152822.
This material is partly based upon work supported by the Swedish Research Council under grant no. 2021-06594 while the third author was in residence at Institut Mittag-Leffler in Djursholm, Sweden during the Fall semester of 2026.}

\begin{document}

\begin{abstract}
We study properties of compact Besicovitch sets and other Kakeya-type sets which are typical in the sense of Baire category.
%
%
%
For classical Besicovitch sets, we prove that a typical compact Besicovitch set $K$ in $\RR^n$  is not sticky, that any line intersects $K$ in a set of ($1$-dimensional) measure zero or in a unit line segment together with a set of measure zero, and that $K$ has no more than $n$ different unit line segments in any direction. Moreover, for each $1\le j\le n$, the set of directions in which $K$ contains exactly $j$ unit line segments has Hausdorff dimension $n-j$.

For compact metric spaces $X$ and $Y$, let $\mathcal{C}(X,Y)$ denote the set of compact subsets of $X\times Y$ whose projection onto the first coordinate is $X$, equipped with the Hausdorff metric.
Kakeya-type sets can be coded by elements of such spaces, allowing typical properties to be studied through comeager subsets of $\mathcal C(X,Y)$. For Besicovitch sets, for example, one can take $X=\mathbb{RP}^{n-1}$ and $Y=[0,1]^n$. We prove a range of properties of typical code sets $A\in \mathcal{C}(X,Y)$ for broad classes of metric spaces $X$ and $Y$, yielding corresponding results for Kakeya-type sets.
Under suitable conditions on $X$ and $Y$, we prove that a typical $A\in \mathcal{C}(X,Y)$ is totally disconnected and purely unrectifiable, its projection to the second coordinate has Hausdorff dimension zero, and the dimension identities $\dim_H(A)=\dim_H(X)$ and $\dim_P(A)=\dim_P(X)+\dim Y$ hold. 

We establish a general sufficient condition under which typical Kakeya-type constructions attain minimal Hausdorff dimension. In fact, under this condition, whenever a construction has zero $s$-dimensional Hausdorff measure for some $s>0$, the same holds for typical constructions.
Thus, in these cases, 
determining the minimal Hausdorff dimension or the existence of a Lebesgue-null construction reduces to studying constructions satisfying any finite or countable collection of typical properties. 
\end{abstract}

\maketitle

\tableofcontents

\section{Introduction}

\subsection{Historical background and motivation}
Kakeya-type problems ask how small a set can be if it contains copies of each element of a given natural collection of sets.
The foundational result in this direction is due to Besicovitch \cite{Besicovitch1919}, who in 1919 constructed a compact set $B$ of Lebesgue measure zero in the plane that contains a unit line segment in every direction. We call a set $E\subseteq \RR^n$ a \emph{Besicovitch set} if it contains a unit line segment in every direction.
The famous Kakeya conjecture states that every Besicovitch set in $\RR^n$ has Hausdorff dimension $n$. This was proved for $n=2$ by Davies \cite{Davies1971} in 1971 and for $n=3$ by Wang and Zahl \cite{WangZahl2025Volume} in 2025.
In 2026, OpenAI \cite{OpenAI2026Kakeya4} released a proof of the $n=4$ case.  Subsequently, OpenAI \cite{OpenAI2026LocalSmoothing} announced a proof of Sogge's local smoothing conjecture for the Euclidean wave equation for dimensions $n\geq2$, from which the Kakeya maximal conjecture, and consequently the Kakeya dimension conjecture, follow in every dimension.

The Kakeya conjecture has motivated a number of related questions. One example is the Besicovitch-$(n,k)$ problem, which asks whether 
every subset of $\RR^n$ containing a $k$-dimensional affine subspace in every direction
must have positive Lebesgue measure, for $2\leq k<n.$
In 2010, Oberlin \cite{Oberlin2010} established this conclusion whenever $n<k+(1+\sqrt{2})^{k-1}$, improving earlier work of Bourgain. More recently, the Kakeya maximal conjecture, together with Oberlin's arguments, implies the positive-measure conclusion for $k>1$. This leaves the study of sets containing $k$-planes in a prescribed and possibly proper subset of directions as a natural related problem.

Other classical Kakeya-type results concern similar questions. In 1976 Stein \cite{Stein1976} proved that, for $n\ge 3$, any subset of $\RR^n$ containing a sphere centered at every point of $\RR^n$ must have positive measure. Ten years later, Bourgain \cite{Bourgain1986} and Marstrand \cite{Marstrand1987} proved the same result for $n=2$. Kakeya-type problems have also been considered over vector spaces other than \(\mathbb{R}^n\). Dvir \cite{Dvir2009} established the finite-field Kakeya conjecture in 2008, while Arsovski \cite{Arsovski2024} resolved its analogue over the \(p\)-adic rationals in 2021. 

In 2003, K\"orner \cite{Korner2003Besicovitch} demonstrated that a Baire category argument can be used to show the existence of small Kakeya-type sets.
He gave a new proof of the existence of zero-measure Besicovitch sets by
showing that, roughly speaking, compact Besicovitch sets of positive measure form a meager subset of the space of all compact Besicovitch sets in the plane.
Thus, a typical compact Besicovitch set has Lebesgue measure zero.
Here and throughout, a property is said to hold typically if the set of exceptions is meager.

Chang, Cs\"ornyei, Héra, and Keleti \cite{Chang2018SmallUnions} 
investigated how small the Hausdorff dimension of a Borel set in $\mathbb R^n$ can be if it contains the $k$-skeleton of a scaled, rotated, or scaled and rotated copy of a cube centered at every point of $\RR^n$.
In all three variants, the upper bounds were obtained by Baire category arguments, showing that typical constructions attain the corresponding lower bounds.
In \cite{Chang2018SmallUnions}, it was also proved that typical $(n,n-1)$-Nikodym sets have Lebesgue measure zero. A similar method established the existence of a Borel set of measure zero in $\RR^n$ containing a hyperplane at every positive distance from every point.
In these results, the copies of the skeletons or the hyperplanes were coded by some compact subsets of a product space. 
The key step was to show that the codes corresponding to sets of positive measure or large Hausdorff dimension form a meager subset of the coding space under the Hausdorff metric topology.
For a more precise description of this framework, see Subsection~\ref{ss:def}.

It was shown in \cite{KeletiMathe2025Kakeya} that typical Besicovitch sets have minimal Hausdorff dimension. 
Consequently, 
any lower bound on Hausdorff dimension established for Besicovitch sets satisfying finitely or countably many typical properties holds for all Besicovitch sets. 
In section~\ref{s:smallest} we extend this minimality principle to a general class of Kakeya-type constructions. More precisely, we introduce a localization property for geometric codings and show that, together with a suitable continuity assumption, it guarantees that typical constructions attain minimal Hausdorff dimension. In fact, we prove the stronger statement that whenever a construction has zero $s$-dimensional Hausdorff measure for some $s>0$, the same holds for typical constructions. We verify the localization property for several classical Kakeya-type codings, including those associated with unit line segments, full lines, and affine planes. This general framework also applies to restricted-direction Kakeya-type problems.

For example, one may ask which compact sets of directions $X\subseteq\mathbf{Gr}_k(\RR^n)$ admit 
a compact family of affine
$k$-planes whose set of directions is $X$ and whose union has
Lebesgue measure zero.
By the results mentioned in the previous paragraph, the existence of a Lebesgue-null construction would imply that typical constructions have zero Lebesgue measure.

Another potential application concerns Mauldin's $P_n$ conjectures, which ask whether every compact purely $1$-unrectifiable set $A\subseteq\RR^2$ with $0<\mathcal H^1(A)<\infty$ intersects some line in at least $n+1$ points. Our results provide a natural class of sets in which to study these questions: typical compact subsets in $[0,1]^2$ whose projection onto the first coordinate is $[0,1]$ are purely $1$-unrectifiable and have $\mathcal H^1(A)=1$. Other possible directions include questions about the minimal Hausdorff dimension of sets containing lower-dimensional spheres centered at every point of a prescribed set, and measure-zero constructions involving translates and dilates of Cantor sets (see \cite{Kel17}, Question~7.1).

As a concrete application, in Section~\ref{s:good} we generalize a planar result of Tam\'as K\'atay \cite{Katay2020} to any dimension by proving the following.

\begin{theorem}
For any $n>1$ there exists a compact collection $\mathcal L$ of lines in $\RR^n$ such that $\mathcal L$ contains a line in every direction and for every line $\Lambda$ in $\R^n$ the intersection $\Lambda \cap \bigcup_{\ell\in \mathcal L, \ell\neq\Lambda} \ell$ is of one-dimensional measure zero.  
\end{theorem}

\subsection{Definitions and Notation}\label{ss:def}

\begin{definition}
Let \(n \geq 2\). If a set \(K \subseteq \mathbb{R}^{n}\) contains a unit line segment in every direction, we call \(K\) a \emph{(classical) Besicovitch set}. If a set \(K \subseteq \mathbb{R}^{n}\) contains a full line in every direction, we call \(K\) a \emph{line Besicovitch set}. If a set \(K \subseteq \mathbb{R}^{n}\) contains a full line of every slope \(m \in [0,1]^{n-1}\), we call \(K\) a \emph{partial (line) Besicovitch set}. Let \(1\leq k < n\). If a set \(K \subseteq \mathbb{R}^n\) contains a \(k\)-dimensional affine subspace in every direction, we call \(K\) an \emph{\((n,k)\)-Besicovitch set}.
\end{definition}

\begin{definition}\label{d:codings}\label{d:gencoding}
    Let $X$ be a compact metric space, $Y$ a separable metric space, and $n\geq 1.$ Define a map $f:X\times Y \to \mathcal P(\mathbb R^n)$ such that $f(x,y)$ is closed in $\mathbb R^n$ for all $(x,y)\in X\times Y$. For $E\subseteq X\times Y$, define 
    \begin{align*}
        K_f(E)= \bigcup_{(x,y)\in E} f(x,y).
    \end{align*}
    We say that $K_f(E)$ is \emph{coded by} $E$.
    We now give four examples of such codings, corresponding to the Besicovitch set types defined above.
    Identify \(\mathbb{RP}^k\) with the space of lines through the origin in $\RR^{k+1}$, and let \(\mathbf{Gr}_j(\mathbb{R}^{n})\) denote the Grassmannian of the \(j\)-dimensional subspaces of \(\mathbb{R}^n\). 
    \begin{enumerate}
    \item \emph{Classical Besicovitch sets.}
    For $v\in\mathbb{RP}^k$, let $l_v$ denote the unit line segment centered at the origin in direction $v$.
    For $E\subseteq\mathbb{RP}^k\times[0,1]^{k+1}$, let
    \begin{align*}
        K_C(E)=\bigcup_{(v,b)\in E}(l_v+b).
    \end{align*}
    \item \emph{Line Besicovitch sets.}
    For $E\subseteq\mathbb{RP}^k\times[0,1]^{k+1}$, let
    \begin{align*}
        K_L(E)=\bigcup_{(L,b)\in E}(L+b).
    \end{align*}
    \item \emph{Partial line Besicovitch sets.}
    For $E\subseteq[0,1]^k\times[0,1]^k$, let
    \begin{align*}
        K_P(E)=\bigcup_{(m,b)\in E} \{(mt+b,t):t\in\RR\}.
    \end{align*}
    \item \emph{$(n,k)$-Besicovitch sets.}
    For $E\subseteq\mathbf{Gr}_k(\RR^n)\times[0,1]^n$, let
\begin{align*}
    K_G(E)=\bigcup_{(L,b)\in E}(L+b).
\end{align*}
\end{enumerate}
In each case, if $E$ has full projection onto its
first-coordinate space, the corresponding coded set is a Besicovitch
set of the corresponding type.
\end{definition}
\begin{remark}
    It is not hard to check that, up to translation and scaling, any compact Besicovitch set contains a modified \(K_C(E)\) for some code set \(E \subseteq  \mathbb{RP}^k \times [0,N]^{k+1}\) with \(N\) sufficiently large. We remark that all our results remain valid for these modified constructions, but for simplicity we state and prove them for \(E \subseteq  \mathbb{RP}^k \times [0,1]^{k+1}\).
    
    Similarly, up to translation and scaling, any line Besicovitch set arising from a compact collection of lines contains a \(K_L(E)\) for some code set \(E \subseteq  \mathbb{RP}^k \times [0,1]^{k+1}\).  
    The analogous statements also hold for partial line Besicovitch sets and \((n,k)\)-Besicovitch sets.
\end{remark}
We now define the spaces of code sets associated with this general framework.
\begin{definition}
\label{1.4}
Let $X$ and $Y$ be metric spaces, and equip \(X \times Y\) with the standard $L_2$ metric for consistency. We define the following spaces of code sets:
\begin{itemize}
    \item If \(X\) is compact and \(Y\) is separable, let 
    \[
    \mathcal{C}(X,Y) \coloneq \{A \subseteq X \times Y: A \text{ is compact and } \mathrm{proj}_1 A = X\}
    \]
    equipped with the Hausdorff metric, where \(\mathrm{proj}_1: X \times Y \to X\) denotes the canonical projection onto the first coordinate.
    \item If \(Y\) is separable, let 
    \[\mathcal{C}_k (Y)\coloneq \{A \subseteq [0,1]^k \times Y: A \text{ is compact and }\mathrm{proj}_1 A = [0,1]^k\}.\]
\end{itemize}
\end{definition}
We consider only the case where \(Y\) is compact, although many of the results can be naturally generalized to the setting where \(Y\) is merely separable. This general formalism is useful as it allows us to study a wide range of Kakeya-type sets such as those discussed in the historical background. 

We now rigorously define what it means for code sets and Kakeya-type sets to be typical.
\begin{definition}
    We say a typical code set \(A \in \mathcal{C} (X,Y)\) satisfies property \(P\) if the subset of \(\mathcal{C}(X,Y)\) satisfying Property \(P\) contains a dense \(G_{\delta}\) set.
\end{definition}
\begin{definition}
\label{1.5}
 For a coding $K_f$ as in Definition $\ref{d:gencoding}$
we call sets of the form $K_f(E)$, where
$E\in\mathcal C(X,Y)$, \emph{$K_f$-Kakeya-type sets}.
We say that a \emph{typical $K_f$-Kakeya-type set
satisfies} property $P$ if a typical
$E\in\mathcal C(X,Y)$ codes a set $K_f(E)$
satisfying $P$.
\end{definition}

As special cases of the above definitions we can talk about properties of \emph{typical Besicovitch sets, line Besicovitch sets, partial Besicovitch sets}, and \emph{\((n,k)\)-Besicovitch sets}.



We conclude with several definitions that are needed to state our main theorems.
\begin{definition}
If \(X\) is a topological space, let \(\dim X\) denote the Lebesgue covering dimension. If \(X\) is a metric space, let \(\dim_H (X)\) and \(\dim_P(X)\) denote its Hausdorff and packing dimensions, respectively.
\end{definition}
\begin{definition}
A topological space \(X\) is \emph{perfect} if it has no isolated points.
\end{definition}
\begin{definition}
    For a linear subspace $L\subseteq\RR^n$, let \(\mathrm{proj}_{L}: \mathbb{R}^n \to L\) denote orthogonal projection onto $L$.
\end{definition}
\begin{definition}
\label{sticky}
    A Besicovitch set \(K \subset \mathbb{R}^n\) is called \emph{sticky} if there exists a set of lines \(L\) of packing dimension \(n-1\), containing at least one line in each direction, such that \(l \cap K\) contains a unit interval for each \(l \in L\). 
\end{definition}

\subsection{Further results}

We now state further results concerning typical
code sets and the four types of Besicovitch
sets introduced above. 
\begin{theorem}
Let \(X\) and \(Y\) be compact metric spaces. A typical \(A\in\mathcal C(X,Y)\) has the following properties:
\begin{enumerate}
    \item (Theorem \ref{2.4}) \(A\) is totally disconnected provided \(Y\) is perfect.
    \item (Corollary \ref{3.6}) \(\max_{x \in X} \# (A \cap (\{x\} \times Y)) = \dim X + 1\) provided \(X\) and \(Y\) are perfect and \(\dim X < \infty\).
    \item (Corollary \ref{4.5}) \(\dim_H (\{x: \# (A \cap (\{x\} \times Y)) = j\}) = k+1-j\) for each \(j \in [k+1]\) provided \(X\) is a Lipschitz manifold with boundary with \(k \coloneq \dim X \geq 1\) and \(Y\) is perfect.
    \item (Proposition \ref{5.2}) For any \(s \in (0,\infty)\), \(\mathcal{H}^s (A) = \mathcal{H}^s (X)\) and \(\dim_H(A) = \dim_H(X)\).
    \item (Proposition \ref{7.3}) \(\dim_H(\mathrm{proj}_2(A)) = 0\).
    \item (Theorem \ref{7main}) \(\dim_P(A) = \dim_P(X) + k\) provided \(\dim_P (X) \in (0,\infty)\) and \(Y = [0,1]^k\).
    \item (Corollary \ref{9.5}) \(A\) is purely unrectifiable provided \(X = [0,1]^k\) and \(Y = [0,1]^n\).
    \item (Theorem \ref{8.3}) If \(X = [0,1]^k\) and \(Y = [0,1]^n\) and \(L\) is a fixed linear subspace of \(\mathrm{span}(X) \times \mathrm{span}(Y)\) with \(L \not\subseteq \mathrm{span}(X) \) then
    \begin{enumerate}
    \item \(\dim_H(\mathrm{proj}_L(A)) = \dim_H(L \cap (\mathrm{span}(X) \times \{0\}))\),
    \item \(\mathrm{proj}_L:A \to L\) is injective.
    \end{enumerate}
    \item (Corollary \ref{10.3}) If \(X = [0,1]^k\) and \(Y = [0,1]^n\) and \(l \subseteq \mathrm{span}(Y)\) is a fixed line through the origin and \(L \coloneq \mathrm{span}(X) \times l\), then \(\mathrm{proj}_L(A)\) is the graph of a continuous surjection \(f: S \to X\) for some \(S \subseteq l\).
\end{enumerate}
\end{theorem}

\begin{theorem}
    A typical Besicovitch set \(K \subset \mathbb{R}^n\) has the following properties:
    \begin{enumerate}
        \item (Corollary \ref{3.13}) \(K\) has no collinear unit line segments.
        \item (Corollary \ref{Cor5.5}) \(K\) has no more than \(n\) different unit line segments in any direction and this upper bound is achieved.
        \item (Corollary \ref{6.6}) Take \(1 \leq j \leq n\). If \(D_j \subseteq \mathbb{RP}^{n-1}\) is the set of directions in which \(K\) has exactly \(j\) unit line segments, then \(\dim_H(D_j) = n-j\).
        \item (Corollary \ref{8.19}) \(K\) is not sticky.
    \end{enumerate}    
\end{theorem}

\begin{theorem}
    A typical line Besicovitch set \(K \subset \mathbb{R}^n\) has the following properties:
    \begin{enumerate}
        \item (Corollary \ref{3.8}) \(K\) is of \(n\)-dimensional Lebesgue measure zero.
        \item (Corollary \ref{3.8}) Let \(K\) be coded by \(A\). For every line \(e \not\in \{l+p: (l,p) \in A\}\), the intersection \(K \cap e\) is of 1-dimensional Lebesgue measure zero.
        \item (Corollary \ref{3.8}) Let \(K\) be coded by \(A\). For every line \(e \in \{l+p: (l,p) \in A\}\), the intersection \(e \cap \bigcup_{(l,p) \in A,l+p \neq e} (l+p)\) is of 1-dimensional Lebesgue measure zero.
        \item (Corollary \ref{line}) \(K\) has no more than \(n\) different lines in any direction and this upper bound is achieved.
        \item (Corollary \ref{6.7}) Take \(1 \leq j \leq n\). If \(D_j \subseteq \mathbb{RP}^{n-1}\) is the set of directions in which \(K\) has exactly \(j\) lines, then \(\dim_H(D_j) = n-j\).
    \end{enumerate}
\end{theorem}

\begin{theorem}
    A typical \((n,k)\)-Besicovitch set \(K \subset \mathbb{R}^n\) has the following properties:
    \begin{enumerate}
        \item (Corollary \ref{5.7}) \(K\) admits a coding such that \(K\) has no more than \(k(n-k)+1\) different coded \(k\)-planes in any direction and this upper bound is achieved.
        \item (Corollary \ref{6.8}) Take \(1 \leq j \leq k(n-k)+1\). The set \(K\) admits a coding such that if \(D_j \subseteq \mathbf{Gr}_k(\mathbb{R}^n)\) is the set of directions in which \(K\) has exactly \(j\) coded \(k\)-planes, then \(\dim_H(D_j) = k(n-k)+1-j\).
    \end{enumerate}
\end{theorem}

\begin{theorem}
    (Corollary \ref{8.5}) Let \(H \subseteq \mathbb{R}^{k+1}\) be any fixed horizontal hyperplane and let \(l \subseteq H\) be an arbitrary line. Then let \(\{h_z\}_{z \in \mathbb{R}} \subseteq H\) be the 1-parameter family of affine \(k-1\)-dimensional hyperplanes of \(H\) perpendicular to \(l\). For a typical partial Besicovitch set \(K \subseteq \mathbb{R}^{k+1}\)
    \begin{enumerate}[label=(\arabic*)]
        \item \(\dim_H(K \cap H) = 0\),
        \item for any \(z \in \mathbb{R}\), \(K\) contains at most one line intersecting \(h_z\).
    \end{enumerate}
\end{theorem}

\subsection{AI declaration}
All the results and proofs were found and written up by the authors. 
ChatGPT-5.6 Sol and ChatGPT-6 Astra helped to find useful references and also to find errors in earlier versions.

\section{Preliminaries}
Throughout the paper, we establish typical properties by combining $G_\delta$ arguments with density constructions. The first three lemmas provide criteria for showing that sets are \(G_{\delta}\), while the latter results develop the covering and approximation tools used to establish density. The first lemma is a folklore result for which we provide a short proof for completeness. Its proof is a simple generalization of the proof of Lemma 2.2 in \cite{Chang2018SmallUnions}.
\begin{definition}
 Let \(\mathcal{K}(X)\) denote the set of nonempty compact subsets of some metric space \(X\) equipped with the Hausdorff metric.
\end{definition}
\begin{lemma}
\label{Lem2.1}Let \(X\) and \(Y\) be arbitrary metric spaces. For any \(s,t \in [0,\infty)\) and continuous map \(F: X \to Y\), the sets
\begin{enumerate}[(i)]
    \item \(S_1 = \{A \in \mathcal{K}(X): \mathcal{H}^s(A) \leq t\}\)
    \item \(S_2 = \{A \in \mathcal{K}(X): \dim_H(A) \leq s \}\)
    \item \(S_3 = \{A \in \mathcal{K}(X): \dim_H(F(A))\leq s \}\)
\end{enumerate}
are \(G_{\delta}\) subsets of \(\mathcal{K}(X)\).
\end{lemma}
\begin{proof}
A set \(A\) has \(s\)-dimensional Hausdorff measure at most \(t\) if and only if for every \(n,m \in \mathbb{N}\), \(A\) has an open \(1/n\)-cover \(\{V_i\}_{i=1}^{\infty}\) such that \(\sum_{i=1}^{\infty} (\diam V_i)^s < t + 1/m\). Since the compact subsets of \(\bigcup_{i=1}^{\infty} V_i\) form an open set, we are done with (i). Then (ii) easily follows. Finally, (iii) follows from (ii) by noting that
\(S_3= F^{-1}(\{A \in \mathcal{K}(Y): \dim_H(A) \leq s\})\). 
\end{proof}
\begin{lemma}
    \label{Lem2.2}Let \(X\) and \(Y\) be compact metric spaces. If for some \(r\in [0,\infty)\),  \(\mathcal{H}^r(X) = s\) and \(\dim_H (X) = t\), then the sets
    \begin{enumerate}[(i)]
    \item \(S_1 = \{A \in \mathcal{C}(X,Y): \mathcal{H}^r(A) = s\}\)
    \item \(S_2 = \{A \in \mathcal{C}(X,Y): \dim_H(A) = t \}\)
    \end{enumerate}
are \(G_{\delta}\) subsets of \(\mathcal{C}(X,Y)\).
\end{lemma}
\begin{proof}
Let \(S_1' \coloneq \{A \in \mathcal{K} (X \times Y): \mathcal{H}^r(A) \leq s\}\). If \(s < \infty\), as \(S_1'\) is \(G_{\delta}\) by Lemma \ref{Lem2.1}, we have that \(S_1' \cap \mathcal{C} (X,Y) = \{A \in \mathcal{C} (X ,Y): \mathcal{H}^r(A) \leq  s\}\) remains \(G_{\delta}\). If \(s = \infty\), then \(S_1' \cap \mathcal{C} (X,Y) = \mathcal{C} (X,Y)\) trivially. Further, as \(\mathrm{proj}_1: X \times Y \to X\) is a 1-Lipschitz map, we have that for all \(A \in \mathcal{C} (X,Y)\), \(\mathcal{H}^r(A) \geq \mathcal{H}^r(X) = s\). Thus \(S_1 = S'_1 \cap \mathcal{C} (X,Y)\) so \(S_1\) is \(G_{\delta}\). 

Analogously, let \(S_2' \coloneq \{A \in \mathcal{K} (X \times Y): \dim_H(A)\leq  t\}\). If \(t < \infty\), as \(S_2'\) is \(G_{\delta}\) by Lemma \ref{Lem2.1}, we have that \(S_2' \cap \mathcal{C} (X,Y)\) is \(G_{\delta}\). If \(t = \infty\), then \(S_2' \cap \mathcal{C}(X,Y) = \mathcal{C}(X,Y)\) trivially. Further, as \(\mathrm{proj}_1: X \times Y \to X\) is Lipschitz, \(\dim_H(A) \geq \dim_H(X) = t\). 
Thus \(S_2 = S_2' \cap \mathcal{C} (X,Y)\), so \(S_2\) is \(G_{\delta}\). 
\end{proof}
\begin{lemma}
    \label{2.3}
    Let \(X\) and \(Y\) be compact metric spaces and let \(\mathcal{D} \subseteq \mathcal{K}(X  \times Y)\) be compact. Fix some integer \(K \geq 2\) and let
    \[G_m(\mathcal{D},K)\coloneq \left\{A \in \mathcal{C} (X,Y): \forall S \in \mathcal{D}, \forall p_1, \dots, p_K \in A \cap S, \min_{i \neq j}\d_{X\times Y} (p_i,p_j) < \frac{1}{m}\right\}.\]
    Then \(G_m\) is an open set.
\end{lemma}
\begin{proof}
    It suffices to show that \(G_m^c\) is closed. Consider some sequence $A_n \in G_m^c$ such that $A_n \rightarrow A$ in the Hausdorff metric. For each $n$, choose \(S_n \in \mathcal{D}\) and \(p_{n,1},\ldots, p_{n,K} \in A_n \cap S_n\)
    such that $d_{X \times Y}(p_{n,i},p_{n,j}) \geq \frac1m$ for all $i\neq j$. By compactness of \(\mathcal{D}\) and \(X \times Y\),
    after passing to a convergent subsequence, \(S_n \to S\) in the Hausdorff metric and
    $p_{n,i} \rightarrow p_i \in S$ for all $i$.
    By Hausdorff convergence, $p_{n,i} \in A_n$ implies $p_i \in A.$ Since limits preserve non-strict inequalities,
    \begin{align*}
    \min_{i \neq j} d_{X \times Y}(p_i,p_j) \ge \frac{1}{m}.
    \end{align*}
    Therefore, $A \in G_m^c$. 
    Thus $G_m^c$ is closed, so $G_m$ is open.
\end{proof}
\begin{definition}
    Given an open cover \(\{A_i\}_{i \in \mathcal{I}}\) of a topological space \(X\), a cover \(\{B_i\}_{i \in \mathcal{I}}\) is a \emph{closed shrinking} if each \(B_i\) is closed in \(X\) and \(B_i \subset A_i\) for all \(i \in \mathcal{I}\).
    \end{definition}
    \begin{definition}
    Let \(\{A_i\}_{i \in \mathcal{I}}\) be a family of subsets of \(X\). 
    An \emph{open swelling} of $\{A_i\}$ is a family of open sets
    \(\{B_i\}_{i \in \mathcal{I}}\) such that \(A_i \subset B_i\) for every \(i \in \mathcal{I}\) and 
    \begin{align*}
        \bigcap_{j=1}^mA_{i_j}\neq\varnothing\iff\bigcap_{j=1}^mB_{i_j}\neq\varnothing
    \end{align*}
    for every finite collection \(i_1,i_2,\dots,i_m \in \mathcal{I}\).
\end{definition}
We have the following useful statements concerning the existence of closed shrinkings and open swellings.
\begin{proposition}
\label{2.2}
\hfill
\begin{enumerate}[a.]
\item \label{a} (Theorem 1.7.8 in Engelking's Dimension Theory \cite{MR1363947}) Every finite open cover \(\{U_i\}_{i=1}^K\) of a normal space \(X\) has a closed shrinking \(\{F_i\}_{i=1}^K\).
\item \label{b} (Theorem 3.1.2 in Engelking's Dimension Theory \cite{MR1363947}) Every finite family \(\{F_i\}_{i=1}^K\) of closed subsets of a normal space \(X\) has an open swelling \(\{V_i\}_{i=1}^K\). Moreover, if there exist some open subsets \(\{U_i\}_{i=1}^K\) of \(X\) satisfying \(F_i \subset U_i\) for all \(i \in [K]\), then the swelling can be defined such that \(\overline{V}_i \subset U_i\). 
\end{enumerate}
\end{proposition}
\begin{lemma}
\label{2.6}
If \(X\) is a perfect normal \(T_1\) space and \(\{U_i\}_{i=1}^M\) is a finite open cover with \(U_j \setminus \bigcup_{i \neq j} U_i \neq \varnothing\) for all \(j \in [M]\), then there exists another open cover \(\{V_i\}_{i=1}^M\) and a closed cover \(\{F_i\}_{i=1}^M\) such that 
$F_j \subset V_j \subset \overline{V}_j \subset U_j$ for all $j\in[M]$,
and 
\begin{align*}
    \# \left( V_j \setminus \bigcup_{i\neq j}V_i\right) = \# \left( F_j \setminus \bigcup_{i\neq j}F_i\right) = \infty.
\end{align*}
\end{lemma}
\begin{proof}
    For each $j\in [M]$, choose \(y_j \in U_j \setminus \bigcup_{i \neq j} U_i\). By Proposition \ref{2.2}\ref{a}, the cover \(\{U_i\}\) has a closed shrinking \(\{F_i\}\). By Proposition \ref{2.2}\ref{b}, \(\{F_i\}_{i=1}^M\) has an open swelling \(\{V_i\}\) satisfying \( F_i\subseteq V_i \subseteq \overline{V}_i \subseteq U_i\) for all \(i \in [M]\). 
    Now, for each $j\in [M]$, define 
    \begin{align*}
        W_j \coloneq V_j \setminus \bigcup_{i\neq j}\overline{V_i}.
    \end{align*}
    Since \(\{V_i\}\) covers \(X\), the point $y_j$ lies in some $V_i$.
    Since $V_i\subseteq U_i$ and $y_j\notin U_i$ if $i\neq j$, $y_j\in V_j.$ This implies that $y_j \notin \overline{V_i}$ when $i\neq j$,
    so $y_j\in W_j.$ Thus $W_j$ is a nonempty open set. Since \(X\) is a perfect \(T_1\) space, this implies that \(\#(W_j) = \infty\). One can easily check that both \(V_j \setminus \bigcup_{i \neq j} V_i\) and \(F_j \setminus \bigcup_{i \neq j} F_i\) contain \(W_j\) so they must also be infinite sets.
\end{proof}
\begin{definition}
    Let \(\{A_i\}_{i\in \mathcal{I}}\) be a cover of a topological space \(X\). The \emph{order} of the cover \(\{A_i\}_{i \in \mathcal{I}}\) is the smallest integer \(k\) such that every \(x \in X\) is contained in at most \(k+1\) covering sets \(A_i\). 
\end{definition}
The following approximation lemma is a key tool for proving density in the space of code sets. It allows us to replace a code set with a finite union of locally constructed pieces while preserving full projection onto $X$ and controlling the Hausdorff distance from the original set.
\begin{lemma}
\label{2.7}
    Let \(X\) and \(Y\) be compact metric spaces, \(E \in \mathcal{C}(X,Y)\) and \(\varepsilon > 0\). Then there exists a finite irredundant closed \(\varepsilon/2\)-cover \(\{F_i\}_{i=1}^N\) of \(X\) and a finite irredundant open 
    \(\varepsilon/2\)-cover \(\{V_j\}_{j=1}^L\) of \(Y\) and an index set \(\mathcal{J} \subseteq [N] \times [L]\) with the following properties
    \begin{enumerate}[(i)]
        \item \(\mathrm{proj}_1(\mathcal{J}) = [N]\).
        \item If we choose nonempty closed subsets \(D(i,j) \subset F_i \times V_j\) for each \((i,j) \in \mathcal{J}\) such that for every \(i \in [N]\) there exists a \(j_i\) with \(\mathrm{proj}_1 (D(i,j_i)) = F_i\), then the set \(E' = \bigcup_{(i,j) \in \mathcal{J}} D(i,j)\) is in \(\mathcal{C} (X,Y)\) and \(\d_H (E, E') \leq \varepsilon\).
    \end{enumerate}
    Furthermore we can also guarantee
    \begin{enumerate}[(a)]
        \item \(\# (V_j \setminus \bigcup_{j' \neq j} V_{j'}) = \infty\) for every \(j \in [L]\), provided \(Y\) is perfect,
        \item \(\# (F_i \setminus \bigcup_{i' \neq i} F_{i'}) = \infty\) for every \(i \in [N]\), provided \(X\) is perfect,
        \item the cover \(\{F_i\}_{i \in [N]}\) has order at most \(\dim X\).
    \end{enumerate}
    Alternatively, 
    \begin{enumerate}[(a)]
        \setcounter{enumi}{3}
        \item for any choice of \(\{V_j\}\) and \(\varepsilon' \in (0,\varepsilon)\) we can choose the cover \(\{F_i\}_{i=1}^N\) such that \(\#\{i' \in [N] \setminus \{i\}: F_i \cap F_{i'}\neq \emptyset\} \leq 3^m-1\) for every \(i \in [N]\) and \(\diam (F_i) \leq \varepsilon'/2\) provided \(X\) is bilipschitz embeddable to \(\mathbb{R}^m\). 
    \end{enumerate}
\end{lemma}
\begin{proof}
    We choose covers \(\{F_i\}_{i=1}^N\) of \(X\) and \(\{V_j\}_{j=1}^L\) of \(Y\) satisfying two key properties: the covers are irredundant and \(\diam (F_i), \diam (V_j) \leq \varepsilon/2\) for all \(i \in [N], j \in [L]\). We also show that statements (a), (b), and (c) or (d) may be guaranteed. By compactness, we may choose finite open \(\varepsilon/4\)-ball covers \(\{B_{\varepsilon/4} (x_i)\}_{i=1}^N\) of \(X\) and \(\{B_{\varepsilon/4} (q_j)\}_{j=1}^L\) of \(Y\) and then let \(\{F_i\}_{i=1}^N \coloneq \{\overline{B}_{\varepsilon/4} (x_i)\}_{i=1}^N\) and \(\{V_j\}_{j=1}^L \coloneq \{B_{\varepsilon/4} (q_j)\}_{j=1}^L\). After discarding redundant members, we may assume these covers are irredundant. Clearly \(\diam (F_i), \diam (V_j) \leq \varepsilon/2\) for all \(i,j\). 

    Next we show (c) by choosing \(\{F_i\}\) differently. If \(\dim X = k\), \(\{B_{\varepsilon/4} (x_i)\}_{i=1}^N\) admits a finite open refinement \(\{U_i\}_{i=1}^K\) with order at most \(k\). By Proposition \ref{2.2}, the open cover \(\{U_i\}_{i=1}^K\) admits a closed shrinking \(\{F_i\}_{i=1}^K\). After deleting redundant covering sets, this cover is irredundant, satisfies \(\diam (F_i) \leq \varepsilon/2\) for all \(i\), and has order at most \(k\).

    Now we show (b) and show that we may obtain (b) and (c) concurrently. As above, if \(\dim X = k\), \(\{B_{\varepsilon/4}(x_i)\}_{i=1}^N\) admits an open refinement \(\{U_i\}_{i=1}^K\) with order at most \(k\). Delete redundant covering sets from \(\{U_i\}_{i=1}^K\). If \(X\) is perfect, by Lemma \ref{2.6}, this refinement admits a closed shrinking \(\{F_i\}_{i=1}^{K}\) satisfying
    \[\# \left(F_i \setminus \bigcup_{i' \neq i} F_{i'} \right) = \infty.\]
    After deleting redundant covering sets, this cover is irredundant, satisfies \(\diam (F_i) \leq \varepsilon/2\) for all \(i\), and satisfies the conditions in statements (b) and (c).

    Next we show (a) by choosing \(\{V_i\}\) differently. If \(Y\) is perfect, by Lemma \ref{2.6}, the cover \(\{B_{\varepsilon/4} (q_j)\}_{j=1}^L\) of \(Y\) induces an open cover \(\{V_j\}_{j=1}^L\) of \(Y\) satisfying
    \[\# \left(V_j \setminus \bigcup_{j' \neq j} V_{j'} \right) = \infty.\]
    After deleting redundant covering sets, this cover is irredundant, satisfies \(\diam (V_j) \leq \varepsilon/2\) for all \(j\), and satisfies the condition in statement (a).

    Next we show (d) by choosing \(\{F_i\}\) differently. After fixing the choice of the cover \(\{V_j\}\), let \(f: X \hookrightarrow \mathbb{R}^m\) be a bilipschitz embedding with inverse Lipschitz constant \(C\). Subdivide \(\mathbb{R}^m\) into a grid of closed axis-parallel cubes \(\mathcal{Q}\) of diameter \(\varepsilon'/(2C)\). As \(X\) is compact, only finitely many \(Q \in \mathcal{Q}\) have \(Q \cap f(X) \neq \emptyset\). Thus let \(\{\tilde{F}_i\}_{i=1}^N\) be the compact cover of \(f(X)\) given by the finitely many nonempty sets of \(\{Q \cap f(X)\}_{Q \in \mathcal{Q}}\). As every \(Q \in \mathcal{Q}\) has nonempty intersection with only its \(3^m-1\) neighbors in \(\mathcal{Q}\), we have that for each \(i \in [N]\)
    \[\#\{i '\in [N] \setminus \{i\}: \tilde{F}_i \cap \tilde{F}_{i'}\neq  \emptyset\} \leq 3^m-1.\]
    Now define \(F_i \coloneq f^{-1}(\tilde{F}_i)\) and obtain a closed cover \(\{F_i\}_{i=1}^N\) of \(X\). After deleting redundant covering sets, this cover is irredundant and as \(f^{-1}\) has Lipschitz constant \(C\), for some \(Q \in \mathcal{Q}\)
    \[\diam F_i \leq C \diam \tilde{F}_i \leq C \diam Q= \varepsilon'/2 < \varepsilon/2\]
    as desired. Further, \(\{F_i\}_{i=1}^N\) inherits the condition in (d) from \(\{\tilde{F}_i\}_{i=1}^N\).

    Now we show statements (i) and (ii). Given covers \(\{F_i\}_{i=1}^N\) and \(\{V_j\}_{j=1}^L\), we define \(\mathcal{J}\) as follows
    \[\mathcal{J} \coloneq \{(i,j) \in [N] \times [L]: (F_i \times V_j) \cap E \neq \emptyset\}.\]
    As \(\mathrm{proj}_1 (E) = X\) and \(\{F_i\}_{i=1}^N\) is an irredundant cover of \(X\), \(\mathrm{proj}_1 (\mathcal{J}) = [N]\). 
    
    Now pick arbitrary nonempty closed subsets \(D(i,j) \subset F_i \times V_j\) for each \((i,j) \in \mathcal{J}\) such that for every \(i \in [N]\) there exists a \(j_i\) with \(\mathrm{proj}_1 D(i,j_i) = F_i\). We then define
    \[E' = \bigcup_{(i,j) \in \mathcal{J}} D(i,j)\]
    as specified by the lemma statement. Clearly \(E'\) is compact, and we further conclude that
    \[\mathrm{proj}_1 (E') \supseteq \mathrm{proj}_1\left(\bigcup_{i =1}^N D(i,j_i)\right) = \bigcup_{i=1}^NF_i = X,\]
    where the \(j_i\) are defined as above. Thus \(E' \in \mathcal{C} (X,Y)\). 

    It remains to show that \(\d_H(E,E') < \varepsilon\). Let \(z \in  E\). Choose \(t \in [N]\), \(s \in [L]\)  such that \(z \in F_t \times V_s\). Then, \((t,s) \in \mathcal{J}\). Now choose \(w \in D(t,s) \subset (F_t \times V_s) \cap E'\). As \(\diam (F_t) \leq \varepsilon/2\), for \(z = (x,y)\), \(w = (x', q)\)
    \[\d_X (x,x') \leq \diam (F_t) \leq \frac{\varepsilon}{2}\]
    and analogously \(\d_Y (y,q) \leq \frac{\varepsilon}{2}\). Under the Euclidean product metric
    \[\d_{X \times Y} (z,w) = \sqrt{(\d_X (x,x'))^2 + (\d_Y (y,q))^2} \leq \frac{\varepsilon}{\sqrt{2}} < \varepsilon.\]
    Hence, \(\d(z,E')<\varepsilon\), so \(E \subseteq B_{\varepsilon} (E')\).

    Conversely let \(w \in E'\). Then \(w \in D(t,s)\) for some \((t,s) \in [N] \times [L]\). By definition of \(\mathcal{J}\), there exists \(z \in E \cap (F_t \times V_s)\). As above this product cell has diameter less than \(\varepsilon\). Since \(z\) and \(w\) lie in the same product cell, \(\d (w, E) < \varepsilon\). Thus we conclude that \(E' \subseteq B_{\varepsilon} (E)\), so \(\d_H (E,E') \leq \varepsilon\).
\end{proof}

\section{Typical Kakeya-type sets are the smallest}\label{s:smallest}
In this section, we establish sufficient conditions under which typical $K_f$-Kakeya-type sets have minimal Hausdorff dimension and, when a measure zero construction exists, have Lebesgue measure zero.

We require two properties of the coding: a continuity property and a localization property.
\begin{definition}
Let $X$, $Y$, $f$ and $K_f$ be as in Definition~\ref{d:gencoding}.

We say that the coding $K_f$ is \emph{continuous} if $f(x,y)$ is compact for every $(x,y)\in X\times Y$ and $f \colon X\times Y \to \cK(\R^n)$ is continuous.

We say that the coding $K_f$ is \emph{$\sigma$-continuous} if there exist continuous functions $f_N\colon X\times Y\to \cK(\R^n)$ ($N\in\NN$) such that for each $(x,y)\in X\times Y$ the sets $f_N(x,y)$ form an increasing chain with union $f(x,y)$.
\end{definition}

It is easy to see that in Definition~\ref{d:codings} the coding of Besicovitch sets is continuous, the others are not continuous (because the sets are not compact) but $\sigma$-continuous.

\begin{notation}
Let $A,B\su \R^n$. 
We write $A\smr B$ if for any $s>0$, $\Hs(B)=0\Rightarrow \Hs(A)=0$.
\end{notation} 

Note that $A\smr B$ implies that $\dim_H(A)\le \dim_H(B)$.
It is also easy to see that this relation is $\sigma$-stable in the sense that $\bigcup_{k=1}^\infty A_k\smr B$ whenever $A_k \smr B$ for every $k$.

\begin{definition}
We say that the coding $K_f$ has the \emph{localization property} if for any $E\in\CXY$ and any nonempty open $V\su Y$ there exists a $D\in\cC(X,V)$ such that $K_f(D)\smr K_f(E)$. 
\end{definition}

\begin{lemma}\label{l:localization}
All four codings of Definition~\ref{d:codings} have the localization property.
\end{lemma}

\begin{proof}
    Fix $E\in \mathcal C(X,Y)$
    and a nonempty open set $V\subseteq Y$. 
    Let us first consider the codings $K_L, K_P,$ and $K_G$.
    In each case, we may apply a sufficiently small dilation
    and a suitable translation to $K_f(E)$ so that all the
    translation parameters lie in $V$. Since dilations and
    translations preserve the directions of lines and planes
    (and the slopes in the case of $K_P$), the resulting
    set is coded by some $D\in\mathcal C(X,V)$.
    Moreover, dilations and translations preserve sets of
    zero Hausdorff measure, so $K_f(D)\smr K_f(E)$.
    The same argument for $K_G$ also applies when the space of directions is any nonempty compact subset of $\mathbf{Gr}_k(\mathbb R^n)$.
    
    For the classical Besicovitch coding $K_C$, let
    $X=\mathbb{RP}^k$ and $Y=[0,1]^{k+1}$. 
    Choose \(q\in V\cap(0,1)^{k+1}\) 
    and $\rho>0$ such that $\overline B_\rho(q)\subseteq V$.
    Cover $Y$ by finitely many compact sets
$B_1,\ldots,B_N$ of diameter less than $\rho$.
Since the closed sets
$\operatorname{proj}_1(E\cap(X\times B_j)),$ $j\in [N]$,
cover $X$, at least one contains a nonempty open
set $U\subseteq X$. Let
$E_0=E\cap(X\times B_j)$ for such a $j$, and fix
$b_0\in\operatorname{proj}_2(E_0)$.

Since the orthogonal group \(O(k+1)\) acts transitively on \(\mathbb{RP}^k\), compactness gives finitely many orthogonal transformations \(R_1,\ldots,R_M\in O(k+1)\) such that
$X=\bigcup_{i=1}^M R_iU$.

For each $i$, apply $R_i$ to $K_C(E_0)$ and
translate the resulting set by $q-R_ib_0$. These rigid motions
preserve the lengths of the coded segments, and
their centers lie in $B_\rho(q)\subseteq V$.
Moreover, the resulting families together contain
a unit segment in every direction. Their union
is therefore coded by some $D\in\mathcal C(X,V)$.

Since this union consists of finitely many
isometric copies of $K_C(E_0)$, and
$K_C(E_0)\subseteq K_C(E)$, we have
$\mathcal H^s(K_C(D)) \leq M\mathcal H^s(K_C(E))$
for every $s>0$. Thus $K_C(D)\smr K_C(E)$,
completing the proof.
\end{proof}
\begin{theorem}\label{l:3.4}
Let $X$ and $Y$ be nonempty compact metric spaces.
Let $K_f$ be a $\sigma$-continuous coding on
$X\times Y$ as in Definition~\ref{d:gencoding}
with the localization property.
Then a typical $E\in\CXY$ has the property that $K_f(E)\smr K_f(A)$ for any $A\in \CXY$.

Consequently, a typical $K_f$-Kakeya-type set has minimal Hausdorff dimension among the $K_f$-Kakeya-type sets.
Furthermore, if there exists a Lebesgue-null $K_f$-Kakeya-type set then a typical $K_f$-Kakeya-type set has Lebesgue measure zero.
\end{theorem}

\begin{proof}
    By $\sigma$-continuity, let
$f_N:X\times Y\to\mathcal K(\mathbb R^n)$,
$N\in\mathbb N$, be continuous maps such that
$f_N(x,y)\subseteq f_{N+1}(x,y)$ and
\[
f(x,y)=\bigcup_{N=1}^\infty f_N(x,y)
\]
for every $(x,y)\in X\times Y$.
Then, for all $s>0$, define $G_s = \{E\in\mathcal C(X,Y):\mathcal H^s(K_f(E)) = 0\}.$
    We first show $G_s$ is a $G_\delta$ set.
    For each $N\in\mathbb N$,
    define
    \begin{align*}
        K_{f,N}(E) = \bigcup_{(x,y)\in E}f_{N}(x,y).
    \end{align*}
    We claim that the map
$\Phi_N:\CXY\to\cK(\R^n)$
defined by $E\mapsto K_{f,N}(E)$
is continuous. Indeed, since $f_N:X\times Y\to\cK(\R^n)$
is continuous, $K_{f,N}(E)$ is compact
when $E$ is compact.
Moreover, for any fixed $E\in\CXY$
and $\varepsilon>0$, continuity of $f_N$
and compactness of $E$ imply that
there exists $\delta>0$ such that, for any $q\in E$,
\begin{align*}
    d_{X\times Y}(p,q)<\delta \implies d_H(f_N(p),f_N(q))<\varepsilon.
\end{align*}
Consequently, if $E'\in\mathcal C(X,Y)$ satisfies $d_H(E,E')<\delta$, then
$d_H(K_{f,N}(E),K_{f,N}(E'))\leq\varepsilon$. Thus $\Phi_N$ is continuous. 

By Lemma $\ref{Lem2.1}$, the set $\{A\in\mathcal K(\mathbb R^n):\mathcal H^s(A)=0\}$
is $G_\delta$. Consequently, 
\begin{align*}
    G_{s,N}=
\{
E\in\mathcal C(X,Y):
\mathcal H^s(K_{f,N}(E))=0
\} 
\end{align*}
is $G_\delta$.
Since the coding is $\sigma$-continuous,
\[
K_f(E)=\bigcup_{N=1}^\infty K_{f,N}(E).
\]
By countable subadditivity, 
\[
G_s=\bigcap_{N=1}^\infty G_{s,N}.
\]
Thus $G_s$ is $G_\delta$.

We now show that if $G_s$ is nonempty, then
$G_s$ is dense. Fix $F\in G_s$,
$E\in\mathcal C(X,Y)$, and $\varepsilon>0$.
Apply Lemma $\ref{2.7}$ to $E$ and $\varepsilon$
to obtain finite covers $\{F_i\}_{i=1}^N$
and $\{V_j\}_{j=1}^L$, together with
an index set $\mathcal J\subseteq[N]\times[L]$.
By the localization property, for every
$j\in[L]$ there exists
$D_j\in\mathcal C(X,V_j)$ such that
$K_f(D_j)\smr K_f(F)$.
For each $(i,j)\in\mathcal J$, let
\[
D(i,j)=D_j\cap(F_i\times Y).
\]
These sets are nonempty and compact,
satisfy $D(i,j)\subseteq F_i\times V_j$,
and have $\operatorname{proj}_1(D(i,j))=F_i$.
Hence, by Lemma~\ref{2.7},
\[
E'=\bigcup_{(i,j)\in\mathcal J}D(i,j)
\]
belongs to $\mathcal C(X,Y)$ and satisfies
$d_H(E,E')\leq\varepsilon$.
Moreover,
\[
K_f(E')\subseteq\bigcup_{j=1}^L K_f(D_j).
\]
Since $F\in G_s$ and
$K_f(D_j)\smr K_f(F)$ for every $j$,
the right-hand side has zero
$s$-dimensional Hausdorff measure.
Thus $E'\in G_s$, proving density.

Let $S=\{s>0:G_s\neq\varnothing\}$.
We claim that
we can choose a countable subset $S_0\subseteq S$
such that
\[
\bigcap_{s\in S}G_s
=
\bigcap_{s\in S_0}G_s.
\]
Indeed, if $S=\varnothing$, take $S_0=\varnothing$.
Otherwise, since $\mathcal H^s(K)=0$ implies
$\mathcal H^t(K)=0$ for every $t>s$, it suffices to choose a sequence in
$S$ converging to $\inf S$, including
$\inf S$ itself if it belongs to $S$.

By the Baire category theorem,
\[
G=\bigcap_{s\in S_0}G_s
\]
is a dense $G_\delta$ subset of $\mathcal C(X,Y)$.
For every $E\in G$ and $A\in\mathcal C(X,Y)$,
if $\mathcal H^s(K_f(A))=0$ for some $s>0$,
then $s\in S$, so $\mathcal H^s(K_f(E))=0$.
Thus $K_f(E)\smr K_f(A)$ for every
$E\in G$ and $A\in\mathcal C(X,Y)$.
\end{proof}
\begin{corollary}
For each coding $K_f\in\{K_C,K_L,K_P,K_G\}$
in Definition~\ref{d:codings}, a typical
$E\in\mathcal C(X,Y)$ satisfies
$K_f(E)\smr K_f(A)$ for every
$A\in\mathcal C(X,Y)$.
Consequently, a typical $K_f$-Kakeya-type set
has minimal Hausdorff dimension among the
$K_f$-Kakeya-type sets, and has Lebesgue
measure zero whenever a Lebesgue-null
$K_f$-Kakeya-type set exists.
\end{corollary}
\begin{proof}
    Immediate from Lemma~\ref{l:localization}, the preceding theorem, and the $\sigma$-continuity of the four codings.
\end{proof}

\section{The Coding Typically has no Bad Redundancies}\label{s:good}
The coding for Besicovitch sets introduced in the introduction has a serious shortcoming. There are two bad ways in which the coding for classical and line Besicovitch sets may be redundant, that is for a code set to code for a Besicovitch set which admits another more \say{natural} coding. 
These redundancies are pathologies of the coding and make proving statements such as Corollaries  \ref{Cor5.5} and \ref{line} impossible. It is therefore essential to show that typical Besicovitch sets do not admit badly redundant codings. There are two types of badly redundant codings.
\begin{enumerate}
    \item For both classical and line Besicovitch sets, it is possible for a compact code set \(A \subseteq \mathbb{RP}^k \times [0,1]^{k+1}\) (resp. \(A \subset [0,1]^k \times [0,1]^{k}\)) to have the property that there exists a closed line segment \(s \subseteq \bigcup_{(l,p) \in A} (l+p)\) (resp. \(s \subseteq K_P(A)\)) such that \(s \not\subset l+p\) for all \((l,p) \in A\) (resp. \(s \not\subset \{(mx+b,x): x \in \mathbb{R}\}\) for all \((m,b) \in A\)).
    \item For line Besicovitch sets, it is possible that for a compact set \(A \subseteq \mathbb{RP}^k \times [0,1]^{k+1}\) there exist distinct \((l,p_1),(l,p_2) \in A\) with \(l + p_1  = l+p_2\). 
\end{enumerate}
We begin by treating bad redundancies of the first type, and in the process, we manage to generalize a two-dimensional result of Tam\'{a}s K\'{a}tay (\cite{Katay2020}) to \(\mathbb{R}^k\).
\begin{definition}
    Let \(\mathcal{L} \) be a compact collection of lines in \(\mathbb{R}^n\). We say that \(\mathcal{L}\) is \emph{good} if for every line \(\Lambda\), the intersection \(\Lambda \cap \bigcup_{l \in \mathcal{L}, l \neq \Lambda } l\) is of one-dimensional measure zero.

    Similarly, let \(A \subseteq \mathbb{RP}^{k} \times [0,1]^{k+1}\). We say that \(A\) is \emph{good} if for every line \(\Lambda\), \(\Lambda \cap \bigcup_{(l,p) \in A, l+p \neq \Lambda}(l+p)\) is of one-dimensional Lebesgue measure zero.

    Finally, let \(A \subseteq [0,1]^k \times [0,1]^{k}\). We say that \(A\) is \emph{good} if for every line \(\Lambda\), \\
    \(\Lambda \cap \bigcup_{(m,b) \in A, \{(mx+b,x): x \in \mathbb{R}\} \neq \Lambda}\{(mx+b,x): x \in \mathbb{R}\}\) is of one-dimensional Lebesgue measure zero.
\end{definition}
\begin{definition}
    We impose the following metric on the space of lines. We let \(\d(l,l') = |p-p'| + \angle (v,v')\), where \(v,v' \in S^{n-1}\) are direction vectors of \(l\) and \(l'\) such that \(\angle (v,v')\) is minimal and \(p\) and \(p'\) are the unique points on \(l\) and \(l'\) such that \(p \perp v\) and \(p' \perp v'\).
\end{definition}
Recall Definition \ref{1.4}.
\begin{proposition}
\label{Katay}
    (\cite{Katay2020}) There exists \(E \in \mathcal{C}_1([0,1])\) which is good.
\end{proposition}
\begin{proposition}
\label{Strongly_good}
For every integer \(n \geq 2\), there exists \(A \in \mathcal{C}_{n-1} ([0,1]^{n-1})\) which is good as well as a compact collection \(\mathcal{L}\) of lines in \(\mathbb{R}^n\) containing a line in every direction which is good.
\end{proposition}
\begin{proof}
By Proposition \ref{Katay}, fix a good \( E \in \mathcal{C}_1([0,1])\). Define
\begin{align*}
  A = \{(m,b): (m_i,b_i) \in E \text{ for every } i \in [n-1]\} \subseteq [0,1]^{n-1} \times [0,1]^{n-1}.   
\end{align*}
We claim that \(A \in \mathcal{C}_{n-1} ([0,1]^{n-1})\) and is good. That \(A \in \mathcal{C}_{n-1} ([0,1]^{n-1})\) is immediate, as it is a Cartesian product of \(n-1\) copies of \(E \in \mathcal{C}_1([0,1])\). 

Now we show that \(A\) is good. Take arbitrary line \(\Lambda\). We consider two cases: the case where \(\Lambda \not\subset H_c \coloneq \{x \in \mathbb{R}^n: x_n = c\}\) for every \(c\in\mathbb R\), and the case where \(\Lambda \subset H_c \) for some \(c\in\mathbb R\). We begin with the first case. 
Then there exist $u,v\in \mathbb R^{n-1}$ such that $\Lambda = \{(ux+v,x):x\in\mathbb R\}$.
Let 
\[X(u,v) \coloneq \bigcup_{(m,b) \in A \setminus \{(u,v)\}}\{(mx+b,x): x \in \mathbb{R}\}.\] 
If $(m,b)\neq (u,v)$ then either $m\neq u$, or $m=u$ and $b\neq v$.
In the latter case, the lines $\{(mx+b, x):x\in\mathbb R\}$ and $\{(ux+v,x):x\in\mathbb R\}$
are disjoint, so it suffices to consider the lines when $m\neq u$.
Then, for at least one $i\in [n-1], m_i\neq u_i.$
For each $i\in[n-1]$, define
\begin{align*}
    T_i\coloneq \{x\in\mathbb R:u_ix+v_i = m_ix+b_i \text{ for some }(m_i,b_i)\in E, m_i\neq u_i\}.
\end{align*}
Since $E$ is good, the line $\{(u_ix+v_i,x):x\in\mathbb R\}$
intersects the union of the other lines coded by $E$ in a set of zero one-dimensional Lebesgue measure.
Since $x\mapsto (u_ix+v_i,x)$ is bilipschitz, $T_i$ has Lebesgue measure zero.
Then 
\begin{align*}
    \{x: (ux+v,x)\in X(u,v)\}\subseteq \bigcup_{i=1}^{n-1}T_i,
\end{align*}
which has Lebesgue measure zero. Thus, $\lambda_1(\Lambda\cap X(u,v))=0.$

Now suppose that \(\Lambda \subset H_c\) for some \(c\in\mathbb R\).
Let \(K_P(E)\) and \(K_P(A)\) be the partial Besicovitch sets coded by \(E\) and \(A\) respectively.
Define
\begin{align*}
    S_c\coloneq \{mc+b:(m,b)\in E\}\subseteq\mathbb R.
\end{align*}
If $h_c\coloneq \mathbb R\times \{c\}\subseteq\mathbb R^2$,
then $h_c\cap K_P(E) = S_c\times\{c\}.$
Since $h_c$ is not coded by $E$ and $E$ is good, $\lambda_1(S_c)=0.$
By definition of $A$, $K_P(A)\cap H_c = S_c^{n-1}\times\{c\}.$
Since $\Lambda$ is a line contained in $H_c$, there exists $i\in [n-1]$ such that 
the coordinate projection $\text{proj}_i:\mathbb R^n\to\mathbb R$
restricts to a bilipschitz map on $\Lambda.$
Further, $\text{proj}_i(\Lambda\cap K_P(A))\subseteq S_c$.
Since $\lambda_1(S_c)=0,\lambda_1(\Lambda\cap K_P(A))=0.$
Thus $A$ is good.

Now let \(\tilde{\mathcal{L}}\) be the compact collection of lines coded by \(A\). Then \(\tilde{\mathcal{L}}\) is good and contains a line of every slope in \([0,1]^{n-1}\). 
By unioning finitely many rotated copies of $\tilde{\mathcal{L}}$,
we may construct a compact collection of lines which is good and contains a line in every direction, as desired.
\end{proof}
Having constructed good collections of lines in every direction, we now show goodness is a typical property. We first need a semicontinuity result for one-dimensional measure under Hausdorff convergence of compact subsets of varying lines.
\begin{lemma}
\label{Lem3.7}
If \(l_i \) is a sequence of lines in \(\mathbb{R}^n\) with \(l_i \to l\) and \(K_i \subset l_i\) is a sequence of compact subsets with one-dimensional Lebesgue measure \(\lambda_1 (K_i) \geq \varepsilon\) with \(K_i \to K \subset l\) compact in the Hausdorff metric, then \(\lambda_1(K) \geq \varepsilon\).
\end{lemma}
\begin{proof}
    For any line \(l\), let \(p(l)\) denote the unique point on \(l\) that is perpendicular to the direction vector of \(l\). Let \(T_i\) be the translation sending \(p(l_i)\) to \(p(l)\). Then let \(\theta_i\) be the minimal angle of rotation needed to rotate \(T_i(l_i)\) to \(l\) about \(p(l)\) and let \(R_i\) be the corresponding rotation about \(p(l)\). Let \(K_i' \coloneq R_i (T_i(K_i)) \subset l\). 
    
    We claim that \(K_i' \to K\) in the Hausdorff metric. As \(l_i \to l\), \(p(l_i) \to p(l)\) and \(\theta_i \to 0\). Hence, \(R_i \circ T_i \to \mathrm{id}\) locally uniformly. As \(K_i \to K\) and \(p(l_i) \to p(l)\), by the triangle inequality, \(\d_H (K_i,p(l))\) remains bounded. Thus all \(K_i \subset \overline{B}_R (p(l))\) for sufficiently large \(R\) and on \(\overline{B}_R(p(l))\), \(R_i T_i \to \mathrm{id}\) uniformly as \(\overline{B}_R(p(l))\) is compact. Thus, we conclude that as \(K_i \to K\), \(K_i' = R_i(T_i(K_i)) \to K\). 
    
    However, each \(K_i'\) is an isometric copy of \(K_i\). Thus, we restrict to \(l\) and conclude by the upper semicontinuity of Lebesgue measure with respect to Hausdorff convergence that \(\lambda_1(K) \geq \varepsilon\) as desired.
\end{proof}
\begin{theorem}
\label{Thm3.8}
A typical \(A \in  \mathcal{C} (\mathbb{RP}^k,[0,1]^{k+1})\) is good. 
\end{theorem}
\begin{proof}
    Let \(K_L(A)\) be the Besicovitch line set coded for \(A \in  \mathcal{C} (\mathbb{RP}^k,[0,1]^{k+1})\). For any line \(l\), let \(p(l)\) denote the unique point on \(l\) that is perpendicular to the direction vector of \(l\). For \(r,\delta,\alpha > 0\), let \(G(r,\delta,\alpha)\) be the set of those \(A \in \mathcal{C} (\mathbb{RP}^k,[0,1]^{k+1})\) for which there does not exist a line \(\Lambda \) such that \(\Lambda \cap \overline{B}_{r}(0) \cap \bigcup_{(l,p) \in A, \angle l, \Lambda \geq \alpha} (l+p)\) has one-dimensional Lebesgue measure at least \(\delta\). Let
    \[G \coloneq \bigcap_{r=1}^{\infty}\bigcap_{m=1}^{\infty}\bigcap_{n=1}^{\infty} G(r, 1/n, 1/m).\]
    Then \(G\) is precisely the set of good sets. It suffices to show that \(G(r,\delta,\alpha)\) is open and that \(G\) is dense.

    We show that \((G(r,\delta,\alpha))^c\) is closed. Consider some sequence \(A_i \in (G(r,\delta,\alpha))^c\) such that \(A_i \to A\) in the Hausdorff metric. For each \(i\), let \(\Lambda_i\) be a line such that \(K_i \coloneq \Lambda_i \cap \overline{B}_{r} \cap \bigcup_{(l,p) \in A_i, \angle l, \Lambda_i \geq \alpha} (l+p)\) has one-dimensional Lebesgue measure at least \(\delta\). It is easy to see that \(K_i\) is compact. As each \(\Lambda_i \cap \overline{B}_{r}(0) \neq \emptyset\) and each \(K_i \subseteq \overline{B}_{r}(0)\), we may extract a subsequence such that \(\Lambda_i \to \Lambda\) and \(K_i \to K \subset \Lambda\) for some line \(\Lambda \) and compact set \(K\). It is easy to see that \(K \subseteq \overline{B}_{r} (0)\). Further, note that \(\bigcup_{(l,p) \in A, \angle l, \Lambda \geq \alpha} (l+p)\) is closed and for every \(\varepsilon > 0\), there exists sufficiently large \(I > 0\) such that for \(i \geq I\),
    \[\overline{B}_{r} (0) \cap \bigcup_{(l,p) \in A_i, \angle l,\Lambda_i \geq \alpha} (l+p)\subseteq B^{\varepsilon}\left( \bigcup_{(l,p) \in A, \angle l, \Lambda \geq \alpha} (l+p)\right).\] 
    We thus conclude that \(K \subseteq \bigcup_{(l,p) \in A, \angle l, \Lambda \geq \alpha} (l+p)\). Therefore
    \[K \subseteq \Lambda \cap \overline{B}_{r}(0) \cap \bigcup_{(l,p) \in A, \angle l, \Lambda \geq \alpha} (l+p).\]

    Now we apply Lemma \ref{Lem3.7} and conclude that as \(\lambda_1(K_i) \geq \delta\), \(\lambda_1(K) \geq \delta\). We conclude that \(A \in (G(r,\delta,\alpha))^c\), so \((G(r,\delta,\alpha))^c\) is closed.

    We now claim that \(G\) is dense. Pick \(E \in \mathcal{C} (\mathbb{RP}^k,[0,1]^{k+1})\) and \(\varepsilon > 0\). Let \(\mathcal{L}\) be the compact collection of lines from Proposition \ref{Strongly_good} which is good. For any line \(l\), let \(L(l)\) be the unique parallel line through the origin. We note that \(S = \{(L(l), p(l)): l \in \mathcal{L}\}\) is compact, so after rescaling and translating the \(y\)-coordinate, we obtain \(S' \in \mathcal{C} (\mathbb{RP}^k,[0,1]^{k+1})\) such that \(\mathcal{L}' \coloneq \{l+p:(l,p) \in S'\}\) is good. 

    Now by Lemma \ref{2.7}, there exists a finite closed \(\varepsilon/2\)-cover \(\{F_i\}_{i=1}^N\) of \(\mathbb{RP}^k\), an open \(\varepsilon/2\)-cover \(\{V_j\}_{j=1}^L\) of \([0,1]^{k+1}\), and an index set \(\mathcal{J} \subseteq [N] \times [L]\) with the following property
    \begin{itemize}
        \item If we choose nonempty closed subsets \(D(i,j) \subset F_i \times V_j\) for each \((i,j) \in \mathcal{J}\) such that for every \(i \in [N]\) there exists a \(j_i\) with \(\mathrm{proj}_1 (D(i,j_i)) = F_i\), then the set \(E' = \bigcup_{(i,j) \in \mathcal{J}} D(i,j)\) is in \(\mathcal{C} (\mathbb{RP}^k,[0,1]^{k+1})\) and \(\d_H (E,E') \leq \varepsilon\).
    \end{itemize}
    For each \((i,j) \in \mathcal{J}\), let \(D(i,j)\) be a copy of \(S' \cap (F_i \times Y)\) with its \(y\)-component scaled and translated such that \(\mathrm{proj}_2D(i,j) \subset V_j\). We obtain that \(E' = \bigcup_{(i,j)} D(i,j)\) is in \(\mathcal{C} (X,Y)\) and \(\d_H (E,E') \leq \varepsilon\). 

    Thus it suffices to show that \(E' \in G\). Pick any line \(\Lambda \). We wish to show that
    \[\lambda_1\left(\Lambda \cap \bigcup_{(l,p) \in E',l+p \neq \Lambda}(l+p)\right) = 0.\]
    As there are only finitely many \(D(i,j)\), it in fact suffices to show that
    \[\lambda_1 \left(\Lambda \cap \bigcup_{(l,p) \in D(i,j),l+p \neq \Lambda} (l+p)\right) =0.\]
    However, \(D(i,j)\) is nothing more than a copy of \(S' \cap (F_i \times Y)\) with the \(y\)-component rescaled and translated. Then let \(\Lambda'\) be the image of \(\Lambda\) under the inverse of the scaling and translating applied above. As \(\mathcal{L}\) is good, we conclude by definition of \(S'\) that
    \begin{align*}
        \lambda_1 \left(\Lambda \cap \bigcup_{(l,p) \in D(i,j),l+p \neq \Lambda} (l+p)\right) & \leq \lambda_1 \left(\Lambda' \cap \bigcup_{(l,p) \in S',l+p \neq \Lambda'} (l+p) \right)
        \\
        &= \lambda_1 \left(\Lambda' \cap \bigcup_{l \in \mathcal{L}', l\neq \Lambda'} l\right) = 0. 
    \end{align*}
    Hence we conclude \(E' \in G\) as desired.
    \end{proof}
    We note that the above theorem is exactly the generalization of Tam\'as K\'atay's theorem in \cite{Katay2020} to higher dimensions, except without the statement that typical line Besicovitch sets are measure zero. This statement has previously been proven for different coding schemes (see \cite{Korner2003Besicovitch}), but for completeness, we note that it follows immediately from Fubini's theorem, the above theorem, and Theorem \ref{3.3}, which we will prove without using this result. We thus have the following full corollary.
    \begin{corollary}
    \label{3.8}
        A typical line Besicovitch set \(B\) in \(\mathbb{R}^n\) coded by a typical code set \(A \in \mathcal{C} (\mathbb{RP}^{n-1},[0,1]^{n})\) satisfies the following three properties
        \begin{enumerate}[(1)]
            \item \(B\) is of \(n\)-dimensional Lebesgue measure zero.
            \item For every line \(e \not\in \{l+p: (l,p) \in A\}\), the intersection \(B \cap e\) is of 1-dimensional Lebesgue measure zero.
            \item For every line \(e \in \{l+p: (l,p) \in A\}\), the intersection \(e \cap \bigcup_{(l,p) \in A,l+p \neq e} (l+p)\) is of 1-dimensional Lebesgue measure zero
        \end{enumerate}
    \end{corollary}

    Now we finish excluding the first class of bad redundancies.
    \begin{proposition}
    \label{Prop3.10}
        A typical \(A \in \mathcal{C} ([0,1]^k,[0,1]^{k})\) is good.
    \end{proposition}
    \begin{proof}
        The proof mostly follows identically to the proof of Theorem \ref{Thm3.8}. The only nontrivial differences are that we instead pick a good  \(S \in \mathcal{C} ([0,1]^k,[0,1]^k)\) by Proposition \ref{Strongly_good} and we instead pick \(\{F_i\}\) to cover \([0,1]^k\). We may then pick \(D(i,j)\) to be a copy of \(S \cap (F_i\times Y)\) with its \(y\)-component rescaled and translated, much as before.
    \end{proof}
    \begin{corollary}
        \label{good1}
        For a typical \(A \in \mathcal{C} (\mathbb{RP}^k , [0,1]^{k+1})\) (resp. \(A \in \mathcal{C}( [0,1]^k, [0,1]^{k})\)), there does not exist a closed line segment \(s \subseteq \bigcup_{(l,p) \in A} (l+p)\) (resp. \(s \subseteq K_P(A)\)) such that \(s \not\subset l+p\) for all \((l,p) \in A\) (resp. \(s \not\subset \{(mx+b,x): x \in \mathbb{R}\}\) for all \((m,b) \in A\)).
    \end{corollary}
    \begin{proof}
        Follows immediately from Theorem \ref{Thm3.8} and Proposition \ref{Prop3.10} by noting that all closed line segments \(s\) have positive one-dimensional Lebesgue measure in the associated line \(\mathrm{aff}(s)\), where \(\mathrm{aff}\) denotes the affine hull. 
    \end{proof}
    We now turn to the second type of redundancy, in which distinct code points represent the same affine subspace. We show this redundancy also disappears for typical code sets, including in the Grassmannian setting.
    \begin{remark}
There are many equivalent choices of metric on the Grassmannian \(\mathbf{Gr}_l(\mathbb{R}^k)\). The one that is most convenient for us is from \cite{barreira2013introduction} and is given by
\[\d (L_1,L_2) \coloneq \d_H (L_1 \cap \overline{B}_1(0), L_2 \cap \overline{B}_1(0)). \]
\end{remark}
\begin{theorem}
\label{4.7}
    Fix \(l \in [k-1]\). For a typical \(A \in \mathcal{C} (\mathbf{Gr}_l(\mathbb{R}^k),[0,1]^{k})\), there do not exist distinct \((L,b),(L,c) \in A\) such that \(L + b = L + c\). 
\end{theorem}
\begin{proof}
    Let 
    \small
    \[G_m \coloneq \left\{ A \in \mathcal{C}(\mathbf{Gr}_l(\mathbb{R}^k),[0,1]^{k}): \forall  L \in \mathbf{Gr}_l(\mathbb{R}^k),(L,b),(L,c) \in A, L+b=L+c \implies |b-c| < \frac{1}{m}\right\}\]
    \normalsize
    and let
    \[G = \bigcap_{m=1}^{\infty} G_m .\]
    Then, the set \(G\) is precisely the set of \(A \in \mathcal{C} (\mathbf{Gr}_l(\mathbb{R}^k), [0,1]^k)\) such that there do not exist distinct \((L,b),(L,c) \in A\) such that \(L + b = L +c\). 
    
    We claim that each \(G_m\) is open. Consider some sequence \(A_n \in G_m^c\) such that \(A_n \to A\) in the Hausdorff metric. For each \(n\), choose \(L_n \in \mathbf{Gr}_l(\mathbb{R}^k)\), \((L_n,b_n),(L_n,c_n) \in A_n \) with \(L_n+b_n = L_n +c_n\) and \(|b_n-c_n|\geq \frac{1}{m}\). By compactness of \(\mathbf{Gr}_l(\mathbb{R}^k)\) and \([0,1]^k\), after passing to a convergent subsequence, we may ensure that \(L_n \to L\), \(b_n \to b\), and \(c_n \to c\) for some \(L \in \mathbf{Gr}_l(\mathbb{R}^k)\) and \(b,c \in [0,1]^k\). By Hausdorff convergence, \((L_n,b_n),(L_n,c_n) \in A_n\) implies \((L,c),(L,b) \in A\). Clearly \(L+b = L+c\) and \(|b-c|\geq \frac{1}{m}\). 
    Therefore, \(A \in G_m^c\). Thus \(G_m^c\) is closed, so \(G_m\) is open.

    Thus it suffices to show that \(G\) is dense. Take \(E \in \mathcal{C} (\mathbf{Gr}_l(\mathbb{R}^k),[0,1]^{k})\) and \(\varepsilon > 0\). We apply Lemma \ref{2.7} to obtain a finite irredundant closed \(\varepsilon/2\)-cover \(\{F_i\}_{i=1}^N\) of \(\mathbf{Gr}_l(\mathbb{R}^k)\), a finite irredundant open \(\varepsilon/2\)-cover \(\{V_j\}_{j=1}^M\) of \([0,1]^{k}\), and an index set \(\mathcal{J} \subseteq [N] \times [M]\) with the following properties
    \begin{enumerate}[(1)]
        \item \(\mathrm{proj}_1 (\mathcal{J}) = [N]\).
        \item If we choose nonempty closed subsets \(D(i,j) \subset F_i \times V_j\) for each \((i,j) \in \mathcal{J}\) such that for every \(i \in [N]\) there exists a \(j_i\) with \(\mathrm{proj}_1 (D(i,j_i)) = F_i\), then the set \(E' = \bigcup_{(i,j) \in \mathcal{J}} D(i,j)\) is in \(\mathcal{C} (X,Y)\) and \(\d_H (E,E') \leq \varepsilon\).
    \end{enumerate}
    Pick \(m\) sufficiently large such that \(\mathbf{Gr}_l(\mathbb{R}^k)\) is bilipschitz embeddable into \(\mathbb{R}^m\) (see Whitney embedding theorem in \cite{GuilleminPollack2010}). Choose \(\varepsilon' \in (0,\varepsilon)\) such that
    \begin{align*}
        (3^m-1)(2\sqrt{k}\varepsilon')^{k-l}k^{l/2}\leq \frac{1}{2} \min_{j \in [M]}\lambda(V_j \cap [0,1]^k),
    \end{align*}
    where \(\lambda\) denotes the Lebesgue measure. By Lemma \ref{2.7}, we can choose the cover \(\{F_i\}\) such that
    \begin{enumerate}[(3)]
        \item \(\#\{i' \in [N] \setminus\{i\}: F_i \cap F_{i'}\neq \emptyset\} \leq 3^m-1\) for every \(i \in [N]\) and \(\diam(F_i) \leq \varepsilon'/2\).
    \end{enumerate}
    
    We now define \(\mathcal{J}_i \coloneq \mathrm{proj}_2(\mathcal{J} \cap (\{i\} \times [M]))\) and pick arbitrary \(j_i \in \mathcal{J}_i\) for each \(i \in [N]\). Our first goal is to select \(y_i \in V_{j_i}\)  for each \(i \in [N]\) inductively such that
    \begin{equation}
        i' < i, L \in F_i \cap F_{i'} \implies L + y_i \neq L+y_{i'}.
        \tag{$\star$}
        \label{star}
    \end{equation}
    We take arbitrary \(i \in [N]\) and suppose we have already selected \(y_{i'}\) for each \(i' \in [i-1]\). We then let \(\mathcal{K}_i \coloneq \{i' \in [N]\setminus \{i\}: F_i \cap F_{i'}\neq \emptyset\}\). By (3), \(|\mathcal{K}_i| \leq 3^m-1\). 
    
    We show that 
    \[S \coloneq (V_{j_i} \cap [0,1]^k)\setminus \bigcup_{i' \in \mathcal{K}_i\cap [i-1]} \bigcup_{L \in F_{i}}(L +y_{i'})\]
    is nonempty by showing \(S\) has positive measure. We fix arbitrary \(L^* \in F_{i}\). As \(\diam (F_i) \leq \varepsilon'\), for every \(L \in F_{i}\), we have that for all \(z \in (L+y_{i'}) \cap [0,1]^{k}\)
    \[\mathrm{dist} (z,L^* +y_{i'})  \leq \d_H(L^*\cap \overline{B}_{\sqrt{k}}(0)  ,L \cap \overline{B}_{\sqrt{k}} (0)) = \sqrt{k} \d(L^*,L)\leq \sqrt{k}\varepsilon'.\]
    This suffices to conclude that
    \begin{align*}
        \bigcup_{L \in F_i} ((L +y_{i'}) \cap [0,1]^k) \subseteq B^{\sqrt{k}\varepsilon'}(L^* +y_{i'}) \cap [0,1]^k,
    \end{align*}
    where \(B^{\sqrt{k}\varepsilon'}(A) \coloneq \bigcup_{x \in A}B_{\sqrt{k}\varepsilon'} (x)\) denotes the generalized ball. 

    Since \(B^{\sqrt{k}\varepsilon'}(L^*+y_{i'}) \cap [0,1]^k\) can be covered by a rotated and translated copy of \([0,2\sqrt{k} \varepsilon']^{k-l} \times [0,\diam ([0,1]^k)]^l\), we have
    \begin{align*}
        \lambda(B^{\sqrt{k}\varepsilon'}(L^* +y_{i'}) \cap [0,1]^k) &\leq \lambda([0,2  \sqrt{k}\varepsilon']^{k-l} \times [0,\diam ([0,1]^k)]^l) = (2\sqrt{k} \varepsilon')^{k-l}k^{l/2}. 
    \end{align*}
    Using \(|\mathcal{K}_i|\leq 3^m-1\) and the choice of \(\varepsilon'\), we conclude that
    \begin{align*}
        \lambda (S) & \geq  \lambda(V_{j_i} \cap [0,1]^k) -|\mathcal{K}_i|(2\sqrt{k}\varepsilon')^{k-l}k^{l/2} \geq \frac{1}{2}\lambda(V_{j_i} \cap [0,1]^k) > 0.
    \end{align*}

    We thus select
    \[y_i \in (V_{j_i}\cap [0,1]^k) \setminus \bigcup_{i' \in \mathcal{K}_i \cap [i-1]} \bigcup_{L  \in F_i}(L +y_{i'}).\]
    We then have that \(L +y_i \neq L + y_{i'}\) for all \(i' \in \mathcal{K}_i \cap [i-1]\) and \(L \in F_i\). Therefore this induction procedure indeed gives \eqref{star}.

    As \(\{F_i\}_{i=1}^N\) is irredundant, we now select arbitrary \(L_i \in F_i \setminus \bigcup_{i' \neq i} F_{i'}\) for each \(i \in [N]\). We next fix \(i \in [N]\) and select \(q_{i,j} \in V_j\) for each \(j \in \mathcal{J}_i \setminus \{j_i\}\) inductively. For each \(j\), we pick
    \[q_{i,j} \in V_{j} \setminus \left(\bigcup_{j' \in (\mathcal{J}_i \setminus \{j_i\}) \cap [j-1]} (L_i + q_{i,j'}) \cup (L_i + y_i)\right). \]
    We now have that \(L_i + q_{i,j} \neq L_i + q_{i,j'}\) for all \( j' \in (\mathcal{J}_i\setminus \{j_i\}) \cap [j-1]\) and \(L_i + q_{i,j} \neq L_i + y_i\).

    We now let \(D(i,j_i) \coloneq F_i \times \{y_i\}\) and let \(D(i,j) \coloneq \{(L_i,q_{i,j})\}\) for \(j \in \mathcal{J}_i \setminus \{j_i\}\). We then let
    \[E' = \bigcup_{(i,j) \in \mathcal{J}} D(i,j).\]
    By (2), \(E' \in \mathcal{C} (X,Y)\) and \(\d_H (E,E') \leq \varepsilon\). Thus it suffices to show \(E' \in G\). Take \(L \in \mathbf{Gr}_l(\mathbb{R}^k)\). Suppose \(L \neq L_i\) for any \(i \in [N]\). Then \(E'
    \cap (\{L\} \times Y) = \{(L,y_i): F_i \ni L \}\). Now \eqref{star} implies that for all distinct \(b_1,b_2 \in \mathrm{proj}_2(E' \cap (\{L\} \times Y))\), we have that \(L + b_1 \neq L+b_2\) as desired.
    
    Now suppose \(L = L_i\) for some \(i \in [N]\). Then \(E' \cap(\{L_i\} \times Y) = \{(L_i,q_{i,j}): j \in \mathcal{J}_i \setminus \{j_i\}\} \cup \{(L_i,y_i)\}\). Once again, by construction, all \(L+q_{i,j} \neq L+q_{i,j'}\) for \(j \neq j'\) and all \(L+q_{i,j} \neq L+y_i\). Thus, \(E' \in G\) as desired. 
\end{proof}
\begin{corollary}
\label{3.13}
    A typical Besicovitch set has no collinear unit line segments.
\end{corollary}
\begin{proof}
Immediate from definitions, Theorem \ref{Thm3.8}, Corollary \ref{good1}, and Theorem \ref{4.7}. 
\end{proof}
\begin{remark}
    The coding for classical Besicovitch sets guarantees distinct points \((v_1,b_1)\), \((v_2,b_2) \in E\) code for distinct line segments \(l_{v_1} + b_1, l_{v_2} + b_2 \subset K_C(E)\). By contrast, the coding for \((k,l)\)-Besicovitch sets (in particular, for line Besicovitch sets) in principle allows \((L,b_1),(L,b_2) \in E\) distinct to code for the same \(l\)-plane \(L + b_1 = L + b_2\). Theorem \ref{4.7} tells us that typically this does not occur. 
\end{remark}

\section{Typical Code Sets are Totally Disconnected}\label{s:disconnected}
The requirement that a code set have full projection onto $X$ does not force it to inherit the connectedness properties of $X$. In fact, we show that typical code sets are totally disconnected whenever $Y$ is perfect.
\begin{theorem}
\label{2.4}
    For compact metric spaces \(X\) and \(Y\) with \(Y\) perfect, a typical \(A \in \mathcal{C} (X,Y)\) is totally disconnected.
\end{theorem}
\begin{proof}
    Define
    \[G_n = \left\{A \in \mathcal{C}(X,Y) : A = \bigsqcup_{i =1}^M A_i, M\in \mathbb{N}, \diam A_i < \frac{1}{n}, A_1, \dots, A_M \text{ compact}\right\},\]
    and let
    \[G = \bigcap_{n=1}^\infty G_n.\]
    If \(A \in G_n\), then each connected component \(K \subset A\) is contained in some \(A_i\), so \(\diam (K) < \frac{1}{n}\). Thus \(G\) is a subset of the set of totally disconnected sets. 
    
    Next we show that \(G_n\) is open. Pick some \(A \in G_n\) and let \(A = \bigsqcup_{i=1}^M A_i\) be the specified decomposition. As the \(A_i\) are compact, we have that \(\mathrm{dist} (A_i,A_j) > 0\) for all \(i \neq j\). For \(M \geq 2\), let
    \[\delta \coloneq  \frac{1}{3} \min \left( \min_{i \neq j} \mathrm{dist} (A_i, A_j), \min_{i} \left( \frac{1}{n} - \diam A_i\right)\right)\]
    and for \(M =1 \), let
    \begin{align*}
        \delta \coloneq \frac{1}{3} \left(\frac{1}{n}-\diam A\right).
    \end{align*}
    One can easily check that \(K \in G_n\) whenever \(\d_H (K,A) < \delta\), so \(G_n\) is open.

    By the Baire category theorem, it suffices to show that each \(G_n\) is dense in \(\mathcal{C} (X,Y)\). Fix \(n \) and consider arbitrary \(E \in \mathcal{C}(X,Y)\) and \(\varepsilon < 1/n\). We construct \(E' \in G_n\) with \(\d_H (E,E') \leq \varepsilon\). 
    
    As \(Y\) is perfect, we apply Lemma \ref{2.7} with statement (a) to obtain a finite closed \(\varepsilon/2\)-cover \(\{F_i\}_{i=1}^N\) of \(X\), a finite open \(\varepsilon/2\)-cover \(\{V_j\}_{j=1}^L\) of \(Y\), and an index set \(\mathcal{J} \subset [N] \times [L]\) with the following properties
    \begin{enumerate}[label=(\arabic*)]
        \item \(\mathrm{proj}_1 (\mathcal{J}) = [N]\).
        \item If we choose nonempty closed subsets \(D(i,j) \subset F_i \times V_j\) for each \((i,j) \in \mathcal{J}\) such that for every \(i \in [N]\) there exists a \(j_i\) with \(\mathrm{proj}_1 (D(i,j_i) )= F_i\), then the set \(E' = \bigcup_{(i,j) \in \mathcal{J}} D(i,j)\) is in \(\mathcal{C} (X,Y)\) and \(\d_H(E,E') \leq \varepsilon\).
        \item \(\# (V_j \setminus \bigcup_{j' \neq j} V_{j'}) = \infty\) for every \(j \in [L]\).
    \end{enumerate}
    As \(\# (V_j \setminus \bigcup_{j' \neq j} V_{j'}) = \infty\), in each \(V_j \setminus \bigcup_{j' \neq j} V_{j'}\) we select \(N\) distinct points \(y_{1,j}, \dots, y_{N,j} \in V_j \setminus \bigcup_{j' \neq j} V_{j'}\). As each \(y_{i,j}\) is contained in only \(V_j\), all \(y_{i,j}\) are distinct. We now take \(D(i,j) \coloneq F_i \times \{y_{i,j}\} \subset F_i \times V_j\) for all \((i,j) \in \mathcal{J}\). If we let \(E' \coloneq \bigcup_{(i,j) \in \mathcal{J}} D(i,j)\), we have that \(E' \in \mathcal{C} (X,Y)\) and \(\d_H(E,E') \leq \varepsilon\) by property (2).

    Thus to show that \(G_n\) is dense, it suffices to show that \(E' \in G_n\). As \(\varepsilon < 1/n\), we have that \(\diam (D(i,j)) = \diam (F_i) \leq \varepsilon < 1/n\). Further, as the \(y_{i,j}\) are distinct, the \(D(i,j)\) are disjoint. Thus \(E' \in G_n\) with the \(D(i,j)\) as the compact disjoint subsets. 
\end{proof}
\section{Vertical Intersections of Code Sets Have Bounded Cardinality}
\begin{theorem}
\label{3.3}
For compact metric spaces \(X\) and \(Y\) with \(X\) perfect and \(\dim X < \infty\), a typical \(A \in \mathcal{C}(X,Y)\) satisfies
\[\# (A \cap (\{x\} \times Y)) \leq \dim X + 1\]
for every \(x \in X\). 
\end{theorem}
\begin{proof}
    Let \(\dim X = k\) and let \(\mathcal{D}= \{\{x\} \times Y: x \in X\}\). Using the notation of Lemma \ref{2.3}, let
    \begin{align*}
    G = \bigcap_{m=1}^\infty G_m(\mathcal{D},k+2).
    \end{align*}
    Then, the set $G $ is precisely the set of all $A \in \mathcal{C}(X,Y)$ whose fibers have cardinality no greater 
    than $k+1$. By Lemma \ref{2.3}, \(G\) is \(G_{\delta}\). Hence, it suffices to prove that $G$ is dense in $\mathcal C(X,Y)$.

    To prove density, fix \(E \in \mathcal{C}(X,Y)\) and \(\varepsilon > 0\). We first treat the case \(k \geq 1\).
    We construct $E'\in G$ with $d_H(E,E')\leq\varepsilon.$ As \(X\) is perfect, we apply Lemma \ref{2.7} with statements (b) and (c) to obtain a finite closed \(\varepsilon/2\)-cover \(\{F_i\}_{i=1}^N\) of \(X\), a finite open \(\varepsilon/2\)-cover \(\{V_j\}_{j=1}^L\) of \(Y\), and an index set \(\mathcal{J} \subseteq [N] \times [L]\) with the following properties
    \begin{enumerate}[label=(\arabic*)]
        \item \(\mathrm{proj}_1 (\mathcal{J}) = [N]\).
        \item If we choose nonempty closed subsets \(D(i,j) \subset F_i \times V_j\) for each \((i,j) \in \mathcal{J}\) such that for every \(i \in [N]\) there exists a \(j_i\) with \(\mathrm{proj}_1 (D(i,j_i) )= F_i\), then the set \(E' = \bigcup_{(i,j) \in \mathcal{J}} D(i,j)\) is in \(\mathcal{C} (X,Y)\) and \(\d_H(E,E') \leq \varepsilon\).
        \item \(\# (F_i \setminus \bigcup_{i' \neq i} F_{i'}) = \infty\) for every \(i \in [N]\).
        \item The cover \(\{F_i\}_{i \in [N]}\) has order at most \(k\).
    \end{enumerate}
    We select \(q_j \in V_j\) for each \(j \in [L]\). As \(\# (F_i \setminus \bigcup_{i' \neq i} F_{i'}) = \infty\), we now select \(L\) distinct points
    \[x_{i,1},\dots,x_{i,L} \in F_i \setminus \bigcup_{i'\neq i} F_{i'}\]
    for each \(i \in [N]\). We define \(H_{i,0} \coloneq F_i\) and \(H_{i,l} \coloneq \{x_{i,l}\}\) for \(1 \leq l \leq L\). This gives a finite closed cover
    \[\{H_{i,l}\}_{i \in [N], l \in [L]_0}\]
    of \(X\), where \([L]_0 \coloneq [L] \cup \{0\}\). Each \(H_{i,l}\) with \(l \neq 0\) intersects no element of the cover besides \(H_{i,0}\). Since \(k \geq 1\), the new cover still has order at most \(k\).

    Let \(\mathcal{J}_i \coloneq \mathcal{J} \cap (\{i\} \times [L])\). As \(\mathcal{J} \subseteq [N] \times [L]\) and \(\mathrm{proj}_1 (\mathcal{J}) = [N]\), for every \(i \in [N]\), \(1 \leq \#(\mathcal{J}_i) \leq L\). Thus, there exists a surjection \(f_i: [L]_0 \to \mathcal{J}_i \) for each \(i \in [N]\). The maps \(f_i\) induce a surjection \(f: [N] \times [L]_0 \to \mathcal{J}\) given by \(f(i,j) = f_i(j)\). We then define
    \[E' \coloneq \bigcup_{i=1}^N \bigcup_{j=0}^L H_{i,j} \times \{q_{\mathrm{proj}_2(f(i,j))}\}.\]
    As \(f_i\) is a surjection  onto \(\mathcal{J}_i\), for every \((i,j) \in \mathcal{J}\) the set \(f^{-1}_i(i,j)\) is nonempty. We thus let
    \[D(i,j) \coloneq \bigcup_{j' \in f_i^{-1}(i,j)} H_{i,j'} \times \{q_{j}\}\]
    for each \((i,j) \in \mathcal{J}\). Then we have
    \[E' = \bigcup_{(i,j) \in \mathcal{J}} D(i,j)\]
    with each \(D(i,j)\) nonempty. For each \(i\), we have that for some \(j_i\), \(f_i(0) = (i,j_i)\), so
    \[D(i,j_i)  \supseteq H_{i,0} \times \{q_{j_i}\}  = F_i \times \{q_{j_i}\}.\]
    Thus for all \(i\), there exists \(j_i\) such that \(\mathrm{proj}_1 (D(i,j_i)) = F_i\), so
    we conclude that \(E' \in \mathcal{C} (X,Y)\) and \(\d_H (E,E') \leq \varepsilon\) by property (2). 

    Thus to show that \(G\) is dense, it suffices to show that \(E' \in G\). We recall our initial definition \(E' = \bigcup_{i=1}^N \bigcup_{j=0}^L H_{i,j} \times \{q_{\mathrm{proj}_2(f(i,j))}\}\). Since the cover \(\{H_{i,j}\}\) has order at most \(k\), all \(x \in X\) lie in at most \(k+1\) different members of \(\{H_{i,l}\}\). Therefore,
    \[\# (E' \cap (\{x\} \times Y)) \leq k+1\]
    for all \(x \in X\), so \(E' \in G\).

    It remains to consider the $k=0$ case.
    Fix $E\in\mathcal C(X,Y)$ and $\varepsilon>0$. As \(X\) is perfect, we apply Lemma \ref{2.7} with statement (c) to obtain a finite disjoint closed \(\varepsilon/2\)-cover \(\{F_i\}_{i=1}^N\) of \(X\) and a finite open \(\varepsilon/2\)-cover \(\{V_j\}_{j=1}^L\) of \(Y\). As the cover \(\{F_i\}\) is disjoint, each $F_i$ is clopen. Now let
    \[\mathcal{J} \coloneq \{(i,j) \in [N] \times [L] : (F_i \times V_j) \cap E \neq \emptyset\}\]
    and let \(\mathcal{J}_i \coloneq \mathcal{J} \cap (\{i\} \times [L]_0)\). Since $\mathrm{proj}_1(E) = X$, each $\mathcal J_i$ is nonempty. Since \(X\) is zero-dimensional and perfect, \(F_i\) can be partitioned into \(\#(\mathcal{J}_i)\) nonempty clopen sets as follows
    \begin{align*}
        F_i = \bigsqcup_{s \in \mathrm{proj}_2(\mathcal J_i)}F_{i,s}.
    \end{align*}
    Pick \(q_j \in V_j\) for each \(V_j\) and define \(D(i,j) \coloneq F_{i,j} \times \{q_j\}\). If we let \(E' \coloneq \bigcup_{(i,j) \in \mathcal{J}} D(i,j)\), we see that \(E'\) is compact and \(\mathrm{proj}_1 (E') = \bigcup_{i=1}^N F_i = X\), so we conclude that \(E' \in \mathcal{C} (X,Y)\). Further, for every \(x \in E\) there exists \(F_i \times V_j \ni x\) with \(x' \in (F_i \times V_j) \cap E'\) and for every \(x' \in E'\) there exists \(F_i \times V_j \ni x'\) with \(x \in (F_i \times V_j) \cap E\). Thus \(\d_H(E,E')  \leq \varepsilon\) follows as \(\diam (F_i \times V_j) \leq \varepsilon\). Finally, as the \(F_{i,j}\) are disjoint, \(E' \in G\), so \(G\) is dense.
\end{proof}
The preceding theorem gives an upper bound on the cardinalities of vertical fibers. To see this bound is sharp, we use a dimension-theoretic result of Hurewicz.
\begin{proposition}
\label{3.4}
    (Hurewicz 1926, \cite{Hurewicz1926}) Let \(X\) be a separable metric space and let \(k\) be a natural number. Then \(\dim X \leq k\) if and only if there exists a zero-dimensional metric space \(Y\) and a continuous and closed surjection \(f: Y \to X\) such that \(\#(f^{-1}(x)) \leq k+1\) for all \(x \in X\).
\end{proposition}
\begin{corollary}
    \label{3.5}
    Let \(X,Y\) be separable metric spaces with \(\dim X = k < \infty\) and let $A \subseteq X \times Y$ be compact and totally disconnected with \(\mathrm{proj}_1 (A) = X\). Then for each \(j \in [k+1]\)
    \[\dim (\{x: \# (A \cap (\{x\} \times Y))\geq j\}) \geq k+1-j\]
    where \(\dim(\varnothing) = -1\) by convention.
\end{corollary}
\begin{proof}
Let
\[S_j \coloneq \{x : \#(A \cap (\{x\} \times Y)) \geq j\}.\]
For the sake of contradiction, assume \(\dim (S_{j})\leq k-j\) for some \(j \in [k+1]\). Then define \(B \coloneq X \setminus S_j\). By the Urysohn-Menger addition formula
\[\dim (B) \geq \dim (X) - \dim (S_j) -1 \geq j-1.\]
We now let \(K = A \cap( B \times Y)\). We then consider the restricted projection
\begin{align*}
    f=\mathrm{proj}_1 |_K&:  K \to B \\
    & (a,b) \mapsto a.
\end{align*}
As \(K \subseteq A\), we have \(\dim (K) = 0\), so we may apply Proposition \ref{3.4} to conclude that there exists \(y \in B\) such that \(|f^{-1} (y)| \geq j\). This, however, contradicts the definition of \(B = X \setminus S_j\), so we have a contradiction. Thus \( \dim (S_j) \geq k+1-j\).
\end{proof}
Combining the lower bound above with the typical upper bound from Theorem $\ref{3.3}$,
we obtain the exact maximal cardinality of a vertical fiber.
\begin{corollary}
    \label{3.6}
    For compact perfect metric spaces \(X\) and \(Y\) with \(\dim X < \infty\), a typical \(A \in \mathcal{C}(X,Y)\) satisfies
    \[\max_{x \in X} \# (A \cap (\{x\} \times Y)) = \dim X + 1.\]
\end{corollary}
\begin{proof}
By Theorem \ref{3.3},  a typical \(A \in \mathcal{C} (X,Y)\) satisfies
\[\max_{x \in X}\# (A \cap (\{x\} \times Y)) \leq \dim X + 1,\]
and  a typical \(A \in \mathcal{C} (X,Y)\) is totally disconnected
by Theorem \ref{2.4}. By taking \(j = k+1\), Corollary \ref{3.5} allows us to conclude that for a typical \(A \in \mathcal{C} (X,Y)\)
\[\max_{x \in X}\# (A \cap (\{x\} \times Y)) = \dim X + 1.\]
\end{proof}
    We now apply these fiber-cardinality bounds to Besicovitch sets. The nonredundancy results of Section~\ref{s:good} let us relate the cardinalities of vertical fibers to the number of line segments in each direction.
\begin{corollary}
\label{Cor5.5}
    A typical Besicovitch set \(K \subseteq \mathbb{R}^k\) has no more than \(k\) different unit line segments in any direction and this upper bound is achieved.
\end{corollary}
\begin{proof}
    Immediate from definitions, Corollary \ref{good1}, Corollary \ref{3.13}, and Corollary \ref{3.6}.
\end{proof}
\begin{corollary}
\label{line}
A typical line Besicovitch set \(K \subseteq \mathbb{R}^{k}\) has no more than \(k\) different lines in any direction and this upper bound is achieved.
\end{corollary}
\begin{proof}
Immediate from definitions, Corollary \ref{good1}, Corollary \ref{3.6}, and Theorem \ref{4.7}.
\end{proof}
\begin{corollary}
\label{5.7}
A typical \((k,l)\)-Besicovitch set \(K \subseteq \mathbb{R}^{k}\) admits a coding such that \(K\) has no more than \(l(k-l)+1\) different coded \(l\)-planes in any direction and this upper bound is achieved.
\end{corollary}
\begin{proof}
Immediate from definitions, Corollary \ref{3.6}, Theorem \ref{4.7}, and the well-known fact that \(\dim(\mathbf{Gr}_l(\mathbb{R}^k)) = l(k-l)\).
\end{proof}
\begin{remark}
Note the stronger statement that a typical \((k,l)\)-Besicovitch set \(K \subseteq \mathbb{R}^{k}\) has no more than \(l(k-l) + 1\) different \emph{arbitrary} \(l\)-planes in any direction does not follow from our results as we lack an analog of Corollary \ref{good1} in the \((k,l)\)-Besicovitch setting. Examining whether analogs of Corollary \ref{good1} or Corollary \ref{line} hold in this setting may be a fruitful direction for future research.
\end{remark}
\section{The Hausdorff Dimensions of Sets of Redundant Directions}
We now refine the cardinality bounds of the previous section by studying the Hausdorff dimension of sets with a given fiber multiplicity. The key tool is a tiling in which 
each face of dimension $k+1-j$ is contained in exactly $j$ tiles, allowing us to obtain the desired dimension bounds.
\begin{lemma}
(\cite{BS20})
    For every integer \(k \geq 1\), there is a face-to-face tiling of \(\mathbb{R}^k\) by congruent convex polytopes of diameter at most 1 such that every \(j\)-dimensional face is contained in exactly \(k+1-j\) full-dimensional tiles.
    \label{4.1}
\end{lemma}
    We first show that the required bounds on the dimensions of fiber-multiplicity sets define a $G_\delta$ condition.
\begin{lemma}
\label{4.2}
    Fix \(k \in \mathbb{N}\). Let \(X\) and \(Y\) be compact metric spaces, let \(F\) be a compact subset of \(X\), and let \(T(F,X,Y) \subseteq \mathcal{C} (X,Y)\) be the set of sets \(A \in \mathcal{C}(X,Y)\) such that for each \(j \in [k+1]\)
    \[\dim_H(\{x \in F: \#(A \cap (\{x\} \times Y)) \geq j\}) \leq k+1-j.\]
    Then \(T(F,X,Y)\) is a \(G_{\delta}\) set.
\end{lemma}
\begin{proof}
    When \(j = 1\), \(\{x \in F: \# (A \cap (\{x\} \times Y)) \geq 1\} = F\) for all \(A \in \mathcal{C} (X,Y)\) due to full projection, so either \(j=1\) provides no constraint or \(T(F,X,Y) = \emptyset\). Fix integer \(j \in [2,k+1]\). For each \(m \in \mathbb{N}\) and \(U \subseteq X\), we define
\begin{align*}
    G(m,U) = \left\{A \in \mathcal{C}(X,Y): \forall x \in F\setminus U,\forall  (x,y_1), \dots, (x,y_{j}) \in A, \min_{i\neq l} \d_Y(y_i,y_l) < \frac{1}{m}\right\}.
\end{align*}
Let \(\mathcal{D} = \{\{x\} \times Y:x \in F \setminus U\}\). If \(U\) is open, \(\mathcal{D}\) is compact, so in the notation of Lemma \ref{2.3}, we have that \(G(m,U) = G_m(\mathcal{D},j)\). Thus, when \(U\) is open, \(G(m,U)\) is open by Lemma \ref{2.3}. We may rephrase the definition of \(G(m,U)\) in more helpful terms. We define
\begin{align*}
    &S_{A,m} \coloneq \left\{x \in F:  \exists (x,y_1), \dots, (x,y_{j}) \in A \text{ distinct}, \min_{i\neq l} \d_Y(y_i,y_l) \geq \frac{1}{m} \right\},\\
    &S_A \coloneq \left\{x \in F:  \exists (x,y_1), \dots , (x,y_{j})\in A\text{ distinct}\right\}.
\end{align*}
The sets \(G(m,U)\) are then given by
\begin{align*}
    G(m,U) = \{A \in \mathcal{C}(X,Y) : S_{A,m} \subseteq U\}.
\end{align*}
Now let
\[U_{p,q} \coloneq \{\text{open } V \subseteq X: \mathcal{H}_{\infty}^{k+1-j+1/p}(V) < 1/q\}.\]
We observe that we may write
\begin{align*}\{\text{compact }A \subseteq X: \dim_H(A) \leq k+1-j\} = \bigcap_{p=1}^{\infty}\bigcap_{q=1}^{\infty}\bigcup_{V \in U_{p,q}}  \mathcal{K}(V).\end{align*}
We replace each \(\mathcal{K}(V)\) with the corresponding \(G(m,V)\) and conclude that
\[\{A \in \mathcal{C}(X,Y): \dim_H (S_{A,m}) \leq k+1-j\} = \bigcap_{p=1}^{\infty}\bigcap_{q=1}^{\infty}\bigcup_{V \in U_{p,q}}  G(m,V). \]
We then intersect over all \(m \in \mathbb{N}\) and conclude that
\begin{align*}
    \{A \in \mathcal{C}(X,Y): \dim_H (S_{A}) \leq k+1-j\} = \bigcap_{m=1}^\infty \bigcap_{p=1}^{\infty}\bigcap_{q=1}^{\infty}\bigcup_{V \in U_{p,q}}  G(m,V).
\end{align*}
Finally, intersecting over all integer \(j \in [2,k+1]\), we conclude that \(T(F,X,Y)\) is a \(G_{\delta}\) set as desired.
\end{proof}
\begin{theorem}
\label{4.3}
    For a compact metric space \(Y\), a typical \(A \in \mathcal{C}_k(Y)\) satisfies the following dimension bound for each \(j \in [k+1]\)
    \[\dim_H (\{x: \#(A \cap (\{x\} \times Y)) \geq j \}) \leq k+1 - j.\]
\end{theorem}
\begin{proof}
    Let \(G\) be the set of sets \(A \in \mathcal{C}_k(Y)\) such that for each \(j \in [k+1]\)
    \[\dim_H(\{x: \#(A \cap (\{x\} \times Y)) \geq j \}) \leq k+1-j.\]
    By Lemma \ref{4.2}, \(G = T([0,1]^k,[0,1]^k,Y)\) is a \(G_{\delta}\) set, so it will suffice to show that \(G\) is dense in \(\mathcal{C}_k(Y)\). Consider arbitrary \(E \in \mathcal{C}_k(Y)\) and \(\varepsilon \in (0,1]\). We will construct \(E'\) with \(\d_H(E,E') \leq \varepsilon\). 
    
    We begin by constructing a closed cover \(\{F_i\}_{i=1}^M\) of \(X\) such that the set of points \(\tilde{I}_j\) contained in \(j\) covering sets has Hausdorff dimension \(k+1-j\) and \(\diam (F_i)\) is sufficiently small. We let \(\{V(\lambda)\}_{\lambda \in \Lambda}\) be the tiling of \(\mathbb{R}^k\) by convex polytopes from Lemma \ref{4.1} where the polytopes are taken to be closed. We let \(I_j\) denote the set of points where \(j\) distinct polytopes intersect. By construction of \(\{V(\lambda)\}_{\lambda \in \Lambda}\), \(I_j\) is contained in the \(k+1-j\)-dimensional faces of the \(\{V(\lambda)\}_{\lambda \in \Lambda}\) so \(\dim_H(I_j) \leq k+1 -j\). 

    As \(Y\) is compact, there exists a finite cover of \(Y\) by balls of radius \(\varepsilon/4\) which we take to be \(\{B_{\varepsilon/4} (q_i)\}_{i=1}^L\). Take \(N \in \mathbb{N}\) sufficiently large such that cubes \(Q\) of side length \(1/N\) have \(\diam Q < \varepsilon/4\). Then decompose the closed unit cube into \(N^k\) small closed cubes of side length \(1/N\) and let \(\mathcal{Q}\) be the collection of these cubes. Choose \(\delta < \varepsilon/4\) sufficiently small such that the interior of each cube \(Q \in \mathcal{Q}\) contains at least \(L\) tiles of \(\{\delta V(\lambda)\}_{\lambda \in \Lambda}\). Let \(F_1, \dots, F_M\) be the nonempty sets of \(\{\delta V(\lambda) \cap [0,1]^k\}_{\lambda \in \Lambda}\). With \(\tilde{I}_j\) defined as above, we have that \(\dim_H(\tilde{I}_j) = \dim_H (I_j)\leq k+1 -j\). 

    For each \(Q \in \mathcal{Q}\), define index sets
    \begin{align*}
        &\mathcal{J}_Q = \{i \in [L]: E \cap (Q \times B_{\varepsilon/4}(q_i)) \neq \varnothing\}, \\
        &\mathcal{K}'_Q = \{i \in [M]: F_i \subseteq \mathrm{int}(Q)\} ,\\
        &\mathcal{K}''_Q = \{i \in [M]: F_i \cap Q \neq \emptyset\}.
    \end{align*}
    As the sets \(\{\mathcal{K}'_Q\}_{Q \in \mathcal{Q}}\) are pairwise disjoint and the sets \(\{\mathcal{K}''_Q\}_{Q \in \mathcal{Q}}\) cover \([M]\), there exist \(\{\mathcal{K}_Q\}_{Q \in \mathcal{Q}}\) with \(\mathcal{K}'_Q \subseteq \mathcal{K}_Q \subseteq \mathcal{K}''_Q\) which partition \([M]\). Then for each \(Q \in \mathcal{Q}\), \(1 \leq |\mathcal{J}_Q|\leq L\leq |\mathcal{K}'_Q|\leq |\mathcal{K}_Q|\), so there exists a surjection \(f_Q: \mathcal{K}_Q \to \mathcal{J}_Q\). Let \(f: [M] \to [L]\) be defined such that \(f|_{\mathcal{K}_Q}= f_Q\). 
    
    We now define our set \(E'\) as
    \[E' = \bigcup_{i=1}^M F_i \times \{q_{f(i)}\}.\]
    As \(E'\) is a finite union of compact sets, \(E'\) is compact. Furthermore, as \(\{F_i\}_{i=1}^M\) covers \([0,1]^k\), \(\mathrm{proj}_1 (E') = [0,1]^k\) so \(E' \in \mathcal{C}_k(Y)\). Now we show that \(E' \in G\). We have that
    \[\{x: \#(E' \cap (\{x\} \times Y)) \geq j\} \subseteq \bigcup_{i=j}^{k+1}\tilde{I}_i,\]
    so we conclude that
    \[\dim_H (\{x: \#(E' \cap (\{x\} \times Y)) \geq j\}) \leq  k+1-j.\]
    Finally \(\d_H (E,E') \leq \varepsilon\) follows analogously as in the proof of Lemma \ref{2.7}, so \(G\) is dense in \(\mathcal{C}_k(Y)\).
\end{proof}
We extend the preceding result to compact Lipschitz manifolds with boundary by applying the construction in local charts, using larger closed neighborhoods to ensure that the unchanged portions of the code set do not interfere with the desired fiber bounds.
\begin{theorem}
\label{4.4}
    For \(X\) a compact Lipschitz manifold with boundary with \(k \coloneq \dim X \geq 1\) and for compact perfect metric space \(Y\), a typical \(A \in \mathcal{C} (X,Y)\) satisfies
    \[\dim_H (\{x: \#(A \cap (\{x\} \times Y)) \geq j\}) \leq k+1-j\]
    for each \(j \in [k+1]\).
\end{theorem} 
\begin{proof}
We begin by constructing a bilipschitz atlas \(\{U_i,\varphi_i\}_{i=1}^M\) with \(F_i\) and \(H_i\) closed such that \(X = \bigcup_{i=1}^M F_i\), \(F_i \subseteq \mathrm{int} (H_i) \subseteq U_i\), and \(\varphi_i(H_i) = [0,1]^k\). As \(X\) is a compact Lipschitz manifold with boundary, \(X\) admits a finite bilipschitz atlas \(\{U_i,\varphi_i\}_{i=1}^M\), that is \(X = \bigcup_{i=1}^M U_i\) and \(\varphi_i: U_i \to (-1,1)^k\) or \(\varphi_i : U_i \to (-1,1)^k \cap \{x_1 \geq 0\} \subseteq \mathbb{R}^k\) are bilipschitz homeomorphisms. By Proposition \ref{2.2}\ref{a}, \(\{U_i\}_{i=1}^M\) has a closed shrinking \(\{F_i\}_{i=1}^M\), so as the \(\varphi_i\) are homeomorphisms, \(\varphi_i (F_i) \) are closed for all \(i \in [M]\). In fact, as \(\varphi_i(F_i)\subseteq (-1,1)^k\), the \(\varphi_i(F_i)\) are compact. We therefore let
\[\delta_i\coloneq \mathrm{dist} (\varphi_i( F_i), \partial (-1,1)^k)/2> 0\]
and we conclude that \(\varphi_i(F_i) \subseteq [-1+2\delta_i,1-2\delta_i]^k\) for each \(i \in [M]\). We let \(r_i \coloneq 1- \delta_i \) and define \(H_i = \varphi_i^{-1} ([-r_i,r_i]^k)\) or \(H_i = \varphi_i^{-1} ([-r_i,r_i]^k \cap \{x_1 \geq 0\})\) which are closed as \(\varphi_i\) is a homeomorphism. As \(\{F_i\}_{i=1}^M\) formed a cover of \(X\), we have that \(\{H_i\}_{i=1}^M\) forms a cover of \(X\). We also have that \(\varphi_i(F_i) \subseteq \mathrm{int} (\varphi_i(H_i))\) (relative to the closed half-space for boundary charts) and \(\varphi_i\) is a homeomorphism so \(F_i \subseteq \mathrm{int}(H_i) \). Finally, as \(\varphi_i\) maps \(H_i\) to closed boxes, we may postcompose with an affine transformation and redefine \(\varphi_i\) such that \(\varphi_i(H_i) =[0,1]^k\) for all \(i \in [M]\). 

Now we let
\[G = \bigcap_{i=1}^M T(F_i,X,Y)\]
where \(T(F_i,X,Y)\) is defined as in the statement of Lemma \ref{4.2}. As the \(F_i\) form a cover of \(X\), \(G\) is precisely the set of sets \(A \in \mathcal{C} (X,Y)\) such that for all \(j \in [k+1]\), \(\dim_H(\{x: \#(A \cap (\{x\} \times Y)) \geq j\}) \leq k+1-j\). It will therefore suffice to show that \(G\) is a dense \(G_{\delta}\) set.

Now each \(T(F_i,X,Y)\) is a \(G_{\delta}\) set by Lemma \ref{4.2}, so by the Baire category theorem it will suffice to show that each \(T(F_i,X,Y)\) is dense. Fix \( i \in [M]\). We seek to push forward the result of Theorem \ref{4.3} to the closed set \(F_i\). However, due to overlaps on the boundary of the chosen chart, we instead perform the construction on \(H_i\) with \(F_i \subseteq \mathrm{int} (H_i)\). 

Take arbitrary \(E \in \mathcal{C}(X,Y)\) and \(\varepsilon > 0\). Define \(\hat{\varphi}_i: H_i \times Y \to [0,1]^k \times Y\) by \(\hat{\varphi}_i (x,y ) = (\varphi_i(x),y)\). We then let
\[\tilde{E}_i \coloneq  \hat{\varphi}_i(E \cap (H_i \times Y)).\]
As \(E \) is compact, \(H_i\) is compact, and \(\hat{\varphi}_i\) is continuous, \(\tilde{E}_i\) is compact. Further, as \(\mathrm{proj}_1 (E) = X\), \(\mathrm{proj}_1(\tilde{E}_i) = \varphi_i(H_i)=[0,1]^k\). Thus \(\tilde{E}_i \in \mathcal{C}_k(Y)\), so we apply Theorem \ref{4.3} to conclude that there exists \(\tilde{E}_i' \in \mathcal{C}_k(Y)\) satisfying 
\[\dim_H(\{z \in [0,1]^k: \#(\tilde{E}'_i \cap (\{z\} \times Y)) \geq j\}) \leq k+1-j\]
for all \(j \in [k+1]\) and \(\d_H (\tilde{E}_i, \tilde{E}_i') < \varepsilon/C\) where \(C\) is the maximum of 1 and the Lipschitz constant of \(\varphi_i^{-1}\). 

We then define \(E'_i \coloneq \hat{\varphi}_i^{-1}(\tilde{E}_i')\) and we finally define
\[E' \coloneq (E \cap ((\mathrm{int} (H_i))^c \times Y) )\cup E'_i.\]
As \(E'\) is the union of compact sets, \(E'\) is compact, and as \(\mathrm{proj}_1 (\tilde{E}_i') = [0,1]^k\), \(\mathrm{proj}_1 (E') = X\). Thus \(E' \in \mathcal{C} (X,Y)\). Further, as \(\d_H (\tilde{E}_i, \tilde{E}'_i)< \varepsilon/C\), \(\d_H (E \cap (H_i \times Y), E_i') < \varepsilon\) as \(\varphi_i^{-1}\) is Lipschitz.

Finally, it is easy to see that for any compact sets \(V, W, Z\) we have that \(\d_H (V \cup Z, W \cup Z)\leq \d_H (V,W)\). Applying this identity for \(V = E \cap (H_i \times Y)\), \(W = E'_i\), and \(Z = E \cap ((\mathrm{int} (H_i))^c \times Y)\) and noting that \(V \cup Z = E\) and \(W \cup Z = E'\), we conclude \(\d_H (E,E') < \varepsilon\). 

We have that the Hausdorff dimension bound for \(E'\) on \(F_i \subseteq \mathrm{int}(H_i)\) is inherited from \(E'_i\), that is for any \(j \in [k+1]\)
\begin{align*}
    \dim_H (\{x \in F_i: \#(E' \cap (\{x\} \times Y)) \geq j\})& \leq \dim_H (\{x \in \mathrm{int}(H_i): \#(E'_i \cap (\{x\} \times Y)) \geq j \}) \\
    & \leq \dim_H (\{z \in [0,1]^k:\#( \tilde{E}_i' \cap (\{z\} \times Y)) \geq j\}) \\
    & \leq k+1-j
\end{align*}
as \(\varphi_i^{-1}\) is Lipschitz. We thus conclude that \(E' \in T(F_i,X,Y)\) so all \(T(F_i,X,Y)\) are dense.
\end{proof}
The upper bounds just established, along with the dimension-theoretic lower bounds from Corollary~\ref{3.5}, determine the Hausdorff dimensions of the sets of points with each possible fiber cardinality.
\begin{corollary}
\label{4.5}
    For \(X\) a compact Lipschitz manifold with boundary with \(k \coloneq \dim X \geq 1\) and for a compact perfect metric space \(Y\), a typical \(A \in \mathcal{C} (X,Y)\) satisfies
    \[\dim_H (\{x: \#(A \cap (\{x\} \times Y)) = j\}) = k+1-j\]
    for each \(j \in [k+1]\).
\end{corollary}
\begin{proof}
Let
\[S_j \coloneq \{x: \#(A \cap (\{x\} \times Y))\geq j\}.\]
A typical \(A \in \mathcal{C}(X,Y)\) is totally disconnected by Theorem \ref{2.4}. Thus, by Corollary \ref{3.5}, \(\dim_H(S_j) \geq \dim (S_j) \geq k+1-j\), so by Theorem \ref{4.4}, we conclude that \(\dim_H(S_j) = k+1-j\). For integer \(j\in [1,k]\)
\begin{align*}k+1-j= \dim_H(S_j) &= \max(\dim_H(S_j \setminus S_{j+1}), \dim_H(S_{j+1})) \\
& = \max(\dim_H (S_j \setminus S_{j+1}), k-j),\end{align*}
so we conclude that
\[\dim_H(\{x: \#(A \cap (\{x\} \times Y)) = j\}) = \dim_H(S_j \setminus S_{j+1}) = k+1-j\]
as desired. When \(j=k+1\), we apply Theorem \ref{3.3} to conclude
\[\dim_H(\{x: \#(A \cap (\{x\} \times Y)) = k+1\}) = \dim_H(S_{k+1}) = 0.\]
\end{proof}
\begin{corollary}
\label{6.6}
    Take \(1 \leq j \leq k\). If \(K \subseteq \mathbb{R}^k\) is a typical Besicovitch set and \(D_j \subseteq \mathbb{RP}^{k-1}\) is the set of directions in which \(K\) has \(j\) unit line segments,  then \(\dim_H(D_j) = k-j\).
\end{corollary}
\begin{proof}
    Immediate from definitions, Corollary \ref{good1}, Corollary \ref{3.13}, and Corollary \ref{4.5}.
\end{proof}
\begin{corollary}
\label{6.7}
Take \(1 \leq j \leq k\). If \(K \subseteq \mathbb{R}^k\) is a typical line Besicovitch set and \(D_j \subseteq \mathbb{RP}^{k-1}\) is the set of directions in which \(K\) has \(j\) lines, then \(\dim_H(D_j) = k-j\).
\end{corollary}
\begin{proof}
    Immediate from definitions, Corollary \ref{good1}, Theorem \ref{4.7}, and Corollary \ref{4.5}.
\end{proof}
\begin{corollary}
\label{6.8}
Take \(1 \leq j \leq l(k-l)+1\). A typical \((k,l)\)-Besicovitch set \(K \subseteq \mathbb{R}^{k}\) admits a coding such that if \(D_j \subseteq \mathbf{Gr}_l(\mathbb{R}^k)\) is the set of directions in which \(K\) has \(j\) coded \(l\)-planes, then \(\dim_H(D_j) = l(k-l)+1-j\).
\end{corollary}
\begin{proof}
    Immediate from definitions, Theorem \ref{4.7}, Corollary \ref{4.5}, and the well-known fact that \(\dim(\mathbf{Gr}_l(\mathbb{R}^k)) = l(k-l)\).
\end{proof}
\section{The Hausdorff Measures of Typical Code Sets}
The projection onto $X$ gives a lower bound for the Hausdorff measure of every code set. To show that this bound is attained typically, we use approximations constructed from closed covers of $X$ whose members have pairwise measure-zero intersection.
\begin{lemma}
\label{CovLem}
    Let \(\mu\) be a nonnegative finite Borel measure on a compact metric space \(X\). For any \(\varepsilon > 0\), there exists a finite closed cover \(\{F_i\}_{i=1}^N\) of \(X\) with \(\diam (F_i) < \varepsilon\) and \(\mu(F_i \cap F_j) = 0\) for all \(i \neq j\). 
\end{lemma}
\begin{proof}
    For every \(x \in X\), consider the family of spheres \(S_r(x)\coloneq \{y \in X: |x-y| = r\}\). As \(X \setminus \{x\}= \bigsqcup_{r \in \mathbb{R}^+} S_r(x)\) and \(\mu(X) < \infty\), there can only be countably many \(r \in \mathbb{R}^+\) such that \(\mu(S_r(x))> 0\). Thus, for each \(x\), pick \(r_x < \varepsilon/2\) such that \(\mu(S_{r_x}(x)) = 0\). Then \(\{B_{r_x} (x)\}_{x \in X}\) is an open cover of \(X\), so by compactness, we may extract an open subcover \(\{B_{r_{x_i}} (x_i)\}_{i=1}^N\). Now define
    \[F_j \coloneq \overline{B}_{r_{x_j}}\setminus \bigcup_{i < j} B_{r_{x_i}}. \]
    Clearly \(\{F_i\}_{i=1}^N\) still covers \(X\). Further, for every \(i \neq j\), we have that \(F_i \cap F_j \subseteq S_{r_{x_i}} (x_i) \cup S_{r_{x_j}} (x_j)\), so \(\mu(F_i \cap F_j) = 0\) as desired. 
\end{proof}
\begin{proposition}
\label{5.2}
    For compact metric spaces \(X\) and \(Y\) and any \(s\in (0,\infty)\), a typical \(A \in \mathcal{C} (X,Y)\) satisfies
    \[\mathcal{H}^s (A) = \mathcal{H}^s(X) \text{ and }\dim_H(A) = \dim_H(X). \]
\end{proposition}
\begin{proof}
Let \(G \coloneq \{A \in \mathcal{C} (X,Y): \mathcal{H}^s(A) =\mathcal{H}^s(X)\}\). By Lemma \ref{Lem2.2}, \(G\) is a \(G_{\delta}\) set.

It thus suffices to show that \(G\) is dense in \(\mathcal{C}(X,Y)\). If \(\mathcal{H}^s(X) = \infty\), then \(G = \mathcal{C} (X,Y)\) as \(\mathrm{proj}_1\) is 1-Lipschitz. If \(\mathcal{H}^s(X) < \infty\), fix \(E \in \mathcal{C} (X,Y)\) and \(\varepsilon > 0\). By compactness, we may choose a finite \(\varepsilon/4\)-ball cover \(\{B_{\varepsilon/4} (q_i)\}_{i=1}^L\) of \(Y\). By applying Lemma \ref{CovLem}, there exists a finite closed cover \(\{F_i\}_{i=1}^N\) of \(X\) with \(\diam (F_i) < \varepsilon/2\) and \(\mathcal{H}^s(F_i \cap F_j) =0\) for all \(i \neq j\). After discarding redundant members, we may assume this cover is irredundant. For each \(t \in [N]\), define the index set
\[\mathcal{J}_t = \{i \in [L]: E \cap (F_t \times B_{\varepsilon/4}(q_i)) \neq \emptyset\}.\]
Since the closed cover of \(X\) was irredundant and \(\mathrm{proj}_1 (E) = X\), each \(\mathcal{J}_t\) is nonempty. Pick a fixed \(j_t \in \mathcal{J}_t\) and an arbitrary \(x_t \in F_t\) for each \(t \in [N]\). Define
\[E' =\bigcup_{t \in [N]} ((F_t \times \{q_{j_t}\}) \cup \{(x_t, q_i): i \in \mathcal{J}_t\}).\]
Clearly \(E'\) is compact and \(\mathrm{proj}_1(E') = X\), so \(E' \in \mathcal{C} (X,Y)\). Further, for every \(x \in E\), there exists \(F_t \times B_{\varepsilon/4} (q_i) \ni x\) with \(x' \in (F_t \times B_{\varepsilon/4} (q_i)) \cap E'\) and for every \(x' \in E'\), there exists \(F_t \times B_{\varepsilon/4} (q_i) \ni x'\) with \(x \in (F_t \times B_{\varepsilon/4} (q_i)) \cap E\). Thus \(\d_H (E,E') < \varepsilon\) follows as \(\diam (F_i \times B_{\varepsilon/4} (q_i))  < \varepsilon\). 

Now we show that \(E' \in G\). As \(s > 0\), the finitely many points \((x_t,q_i)\) contribute no measure. Using that \(\mathrm{proj}_1\) is 1-Lipschitz and \(\mathcal{H}^s (F_i \cap F_j) = 0\) for all \(i \neq j\), we obtain
\[\mathcal{H}^s(E') = \mathcal{H}^s\left(\bigcup_{t \in [N]} (F_t \times \{q_{j_t}\} )\right) = \sum_{t \in [N]} \mathcal{H}^s(F_t) = \mathcal{H}^s\left(\bigcup_{t \in [N]} F_t\right) \leq \mathcal{H}^s (X),\]
so \(E' \in G\).

That a typical \(A \in \mathcal{C}(X,Y)\) satisfies \(\dim_H(A) = \dim_H(X)\) follows from the fact that
\[\{A \in \mathcal{C} (X,Y): \dim_H(A) = \dim_H(X)\} \supseteq \bigcap_{s \in \mathbb{Q}^+} \{A \in \mathcal{C}(X,Y): \mathcal{H}^s (A) = \mathcal{H}^s (X)\}\]
contains a dense \(G_{\delta}\) set by the Baire category theorem.
\end{proof}
In contrast to the full first-coordinate projection, the second-coordinate projection of a typical code set can be much smaller. We show that it has Hausdorff dimension zero.
\begin{proposition}
\label{7.3}
    For compact metric spaces \(X\) and \(Y\), a typical \(A \in \mathcal{C} (X,Y)\) satisfies
    \[\dim_H(\mathrm{proj}_2(A)) = 0.\]
\end{proposition}
\begin{proof}
Let \(G \coloneq  \{A \in \mathcal{C} (X,Y): \dim_H(\mathrm{proj}_2(A)) = 0\}\).
By Lemma \ref{Lem2.1}, \(G\) is a \(G_{\delta}\) set. Thus it suffices to show that \(G\) is dense.

Fix \(E \in \mathcal{C} (X,Y)\) and \(\varepsilon > 0\). By Lemma \ref{2.7}, we obtain a finite closed \(\varepsilon/2\)-cover \(\{F_i\}_{i=1}^N\) of \(X\), a finite open \(\varepsilon/2\)-cover \(\{V_j\}_{j=1}^L\) of \(Y\), and an index set \(\mathcal{J} \subset [N] \times [L]\) with the following property
\begin{itemize}
        \item If we choose nonempty closed subsets \(D(i,j) \subset F_i \times V_j\) for each \((i,j) \in \mathcal{J}\) such that for every \(i \in [N]\) there exists a \(j_i\) with \(\mathrm{proj}_1 (D(i,j_i) )= F_i\), then the set \(E' = \bigcup_{(i,j) \in \mathcal{J}} D(i,j)\) is in \(\mathcal{C} (X,Y)\) and \(\d_H(E,E') \leq\varepsilon\).
\end{itemize}
We choose \(q_j \in V_j\) for each \(j \in [L]\) and let \(D(i,j) \coloneq F_i \times \{q_j\}\) for each \(i,j\). We clearly have that all \(D(i,j)\) are nonempty and all \(D(i,j)\) satisfy \(\mathrm{proj}_1 (D(i,j)) = F_i\). Thus, if we let \(E' \coloneq \bigcup_{(i,j) \in \mathcal{J}} D(i,j)\), we immediately conclude that \(E' \in \mathcal{C} (X,Y)\) and \(\d_H(E,E') \leq \varepsilon\).

It therefore suffices to show that \(E' \in G\). However, \(\mathrm{proj}_2(E') \subseteq \{q_j:j \in [L]\}\) which is a finite collection of points, so \(\dim_H(\mathrm{proj}_2(E')) = 0\) as desired.
\end{proof}
\section{Typical Besicovitch Sets are not Sticky}
Recall Definition \ref{sticky}. The distinction between sticky and non-sticky Besicovitch sets was crucially important for Wang and Zahl's proof of the 3D Kakeya conjecture (\cite{WangZahl2025Volume}). Here we give a proof of the result that typical Besicovitch sets are not sticky, and in doing so show that for typical \(A \in \mathcal{C} (X,[0,1]^k)\), \(\dim_P(A) = \dim_P(X) + k\). 
To transfer typical properties from code sets over Cantor subsets of $X$ to code sets over $X$ itself, we first establish that restriction preserves residual properties.
\begin{definition}
    Let \(X\) and \(Y\) be topological spaces. A map \(f: X \to Y\) is \emph{quasicontinuous} if for every \(x \in X\) and open sets \(U \ni x\) and \(V \ni f(x)\), there exists some nonempty open \(W \subseteq U\) with \(f(W) \subseteq V\).
\end{definition}
\begin{proposition}
\label{quopen}(\cite{Crannell2006})
    Let \(X\) and \(Y\) be compact metric spaces and let \(f: X \to Y\) be open and quasicontinuous. If \(R \subset Y\) is residual in \(Y\), then \(f^{-1} (R) \subset X\) is residual in \(X\). 
\end{proposition}
\begin{lemma}
\label{restriction}
    Let \(X\) and \(Y\) be compact metric spaces and let \(H \subseteq X\) be a Cantor set. Let \(R: \mathcal{C} (X,Y) \to \mathcal{C} (H,Y)\) be the restriction map \(R(A) = A \cap (H \times Y)\). If \(G \subset \mathcal{C}(H,Y)\) is residual in \(\mathcal{C} (H,Y)\), then \(R^{-1} (G)\) is residual in \(\mathcal{C} (X,Y)\).
\end{lemma}
\begin{proof}
By Proposition \ref{quopen}, it suffices to show that \(R\) is open and quasicontinuous. To show that \(R\) is open, it suffices to show that for every \(A \in \mathcal{C} (X,Y)\) and \(\varepsilon > 0\), there exists \(\delta > 0\) such that \(B_{\delta} (R(A)) \subseteq R(B_{\varepsilon} (A))\). 

Fix \(A \in \mathcal{C} (X,Y)\) and \(\varepsilon > 0\) and let \(\delta \coloneq \varepsilon/3\). Take some \(D\in B_{\delta} (R(A))\). We construct \(D' \in B_{\varepsilon}(A)\) with \(D= R(D')\). Let \(G\) be the union of those open balls \(B_r(p)\) of \(X \times Y\) which are disjoint with \(D\) and for which \(p\in R(A)\) and \(r<\delta\). Clearly \(G\) is open and disjoint from \(D\). As \(D\) is closed, \(R((A \cup D) \setminus G) = D\). Further, \((A \cup D) \setminus G\) is trivially a compact subset of the generalized ball \(B^{\delta}(A)\). Also, \(G \subseteq B^{\delta}_{X \times Y}(R(A)) \subseteq B^{\varepsilon}_{X \times Y}(D)\), so \(A \subseteq B^{\varepsilon}((A \cup D) \setminus G)\). The set \((A \cup D) \setminus G \) almost suffices as our choice of \(D'\), but it does not have full projection. We therefore add an auxiliary set 
\[K \coloneq \{z \in X \times Y: \mathrm{dist} (z,D) = \mathrm{dist} (z,H \times Y) \leq \delta \}. \]
Clearly \(K\) is compact. We claim that \(\mathrm{proj}_1(K) =\overline{B}^{\delta} (H)\). Pick \(x \in \overline{B}^{\delta} (H)\). As \(H\) is compact, there exists a minimizer \(h \in H\) such that \(\d(x,h) = \d (x,H) \leq \delta\). As \(\mathrm{proj}_1 (D) = H\), there exists \(y \in Y\) with \((h,y) \in D\). Clearly \(\mathrm{dist} ((x,y),D)= \mathrm{dist} ((x,y),H \times Y) =\d (x,h)\), so there exists \((x,y) \in K\) for each \(x \in \overline{B}^{\delta} (H)\).

By definition of \(K\), we clearly have that \(K\subseteq \overline{B}^{\delta}_{X \times Y} (D) \subseteq B^{\varepsilon} (A)\). Now define \(D'\coloneq K \cup (A \cup D) \setminus G \). We observe that \(D' \in B_{\varepsilon} (A)\) and \(R(D') = D\) as desired. 

Now we show that \(R\) is quasicontinuous. By Theorem \ref{3.3}, as \(H\) is a Cantor set, the set of \(A \in \mathcal{C}(H,Y)\) with \(\#(A \cap (\{x\} \times Y)) = 1\) for all \(x \in H\) is dense in \(\mathcal{C}(H,Y)\). If a set \(A\) has this property, then \(A\) is a compact graph of some function \(f: H \to Y\). Thus, this function \(f\) is continuous. Further, Lipschitz functions are dense in the space of continuous functions over a Cantor set, so we conclude that the set 
\[\mathcal{L}\coloneq \{A \in \mathcal{C}(H,Y): A \text{ is the graph of a Lipschitz function }f: H \to Y\}\] 
is dense in \(\mathcal{C}(H,Y)\).

Let \(R: \mathcal{C}(X,Y) \to \mathcal{C}(H,Y)\) be given by \(R(A) = A \cap (H \times Y)\) as before. Pick any \(A \in \mathcal{C}(X,Y)\) and suppose that \(U \ni A\) and \(V \ni R(A)\) are two open sets. To show that \(R\) is quasicontinuous, it suffices to show that there exists open nonempty \(W \subseteq U\) with \(R(W) \subseteq V\). As \(R\) is open, \(V \cap R(U)\) is open in \(\mathcal{C}(H,Y)\). Hence, as \(\mathcal{L}\) is dense, there exists \(K \in \mathcal{L}\) such that \(K \in V \cap R(U)\). As \(K \in R(U)\), there exists \(K' \in U\) with \(K = R(K')\). Let \(f\) be the Lipschitz function corresponding to \(K\).

We now claim that for every \(\varepsilon > 0\) there exists \(\delta > 0\) such that \(K' \cap (B^{\delta}(H) \times Y) \subseteq B^{\varepsilon}_{X \times Y}(K)\). Suppose the contrary. Then for every \(n \in \mathbb{N}\), we may find \(x_n \in (K' \cap (B^{1/n} (H) \times Y)) \cap (B^{\varepsilon}(K))^c\). By compactness of \(K'\), we extract a convergent subsequence with \(x_{k_n} \to x \in K'\). As \(H\) is closed \(x \in K'\cap (H \times Y)\), but also \(x \in K^c\), so we have a contradiction.

Now let \(\varepsilon > 0\) be sufficiently small such that \(B_{\varepsilon}(K) \subseteq V\). Let \(C\) be the Lipschitz constant of \(f\) and let \(\tilde{C} = \max(1,C)\). Pick \(\delta \in (0,\varepsilon/(10\tilde{C}))\) sufficiently small such that \(B_{\delta} (K') \subseteq U\) and \(K' \cap (B^{\delta}(H) \times Y) \subseteq B^{\varepsilon/(10\tilde{C})}_{X \times Y}(K)\). Suppose \(E \in B_{\delta}(K')\). For each \(x \in H\), pick \(y \in \mathrm{proj}_2(E \cap ({x} \times Y))\). Then there must exist some \(z \in K' \cap (B_{\delta}(x) \times Y)\) with \(\d((x,y),z) < \delta\). However, we now have that \(z \in K' \cap (B^{\delta}(H) \times Y) \subseteq B^{\varepsilon/(10\tilde{C})}_{X \times Y}(K)\). Specifically, as \(z \in K' \cap (B^{\delta}(x) \times Y)\), we must have that there exists \(x' \in H \cap B_{\varepsilon/(5\tilde{C})}(x)\) such that \(z \in B_{\varepsilon/(10\tilde{C})}(x',f(x'))\), so \((x,y)\in B_{\varepsilon/(5\tilde{C})}(x',f(x'))\). However, as \(f\) is Lipschitz with constant \(C\), we conclude that \((x,y) \in B_{\varepsilon} (x,f(x))\). 

Applying this argument for each \(x \in H\) and \(y \in \mathrm{proj}_2(E \cap (\{x\} \times Y))\) suffices to show that \(R(E) \in B_{\varepsilon}(K) \subseteq V\), so we may choose \(W = B_{\delta} (K')\). Thus \(R\) is quasicontinuous as desired.
\end{proof}
\begin{definition}
    If \(X\) is a metric space, we let \(\overline{\dim}_{B} (X)\) denote the upper box dimension and we let \(\underline{\dim}_{B}(X)\) denote the lower box dimension. 
\end{definition}
\begin{definition}
    If \(X\) is a complete metric space and \(A \subseteq X\), the \emph{smallest local box dimension} of the set \(A\) is given by
    \[\mathrm{sl-}\overline{\dim}_B (A) \coloneq \inf \{\overline{\dim}_B (\overline{B}_r(x) \cap A): x \in A, r> 0\}. \]
\end{definition}
\begin{lemma}
\label{Lem7.3}
(\cite{MyjakRudnicki2000}) If \(X\) is a complete metric space and \(A \subseteq X\) is closed, then
\(\dim_{P}(A)\geq \mathrm{sl-}\overline{\dim}_B(A)\). 
\end{lemma}
\begin{lemma}
\label{Lem7.4}
Let \(X\) be a separable metric space. Then for any \(s > 0\) the set
\[S = \{A \in \mathcal{K}(X) : \mathrm{sl-}\overline{\dim}_B (A) \ge s\}\]
is a \(G_{\delta}\) subset of \(\mathcal{K}(X)\).
\end{lemma}
\begin{proof}
As \(X\) is separable, \(X\) has a countable dense subset \(D \subseteq X\). Let \(A \subseteq X\). For every \(x \in A, r \in (0,\infty)\), we may find \(y \in D, q \in \mathbb{Q}^+\) with \(B_q(y) \subseteq \overline{B}_r(x)\) and for every \(y \in D, q \in \mathbb{Q}^+\) with \(B_q(y) \cap A \neq \emptyset\), we may find \(x \in A, r \in \mathbb{R}\) with \(\overline{B}_r(x) \subseteq B_q(y)\). Thus
\[\mathrm{sl-}\overline{\dim}_B(A ) = \inf \{\overline{\dim}_B(B_q(y) \cap A): y \in D, q \in \mathbb{Q}^+, B_q(y) \cap A \neq \emptyset\}.\]
Let \(N(A,\varepsilon)\) denote the maximum number of disjoint closed balls of radius \(\varepsilon\) that may be found with centers in \(A\). We claim that for any \(q \in \mathbb{Q}^+\), \(y \in D\), \(\varepsilon > 0\), and \(k > 1/s\) the set
\[G(q,y,\varepsilon,k) = \left\{A \in \mathcal{K}(X):\exists \varepsilon' \in (0,\varepsilon] \text{ such that }N(B_q(y) \cap A , \varepsilon') > (\varepsilon')^{-\left(s- \frac{1}{k}\right)}\right\}\]
is open. 

Suppose \(A \in G(q,y,\varepsilon,k)\). Then there exists \(\varepsilon' \in (0,\varepsilon]\) and a collection of disjoint balls \(\{\overline{B}_{\varepsilon'} (x_i)\}_{i=1}^N\) with \(x_i \in B_q(y) \cap A\) and \(N > (\varepsilon')^{-\left(s- \frac{1}{k}\right)}\). For each \(x_i\) there exists \(\delta_i> 0\) sufficiently small such that \(B_{\delta_i} (x_i) \subset B_q(y)\). Let \(\delta = \min(\varepsilon',\min_i \delta_i)\). Choose \(\beta \in (0,\delta)\) small enough such that \(N > (\varepsilon'-\beta)^{-\left(s-\frac{1}{k}\right)} \). Take \(A' \in B_{\beta} (A)\). For each \(x_i\), pick \(x'_i \in B_{\beta} (x_i) \cap A'\). As \(\beta< \delta\), all \(x'_i \in B_q(y)\). We then have that \(\overline{B}_{\varepsilon' - \beta} (x'_i) \subseteq \overline{B}_{\varepsilon'} (x_i)\), so the balls \(\{\overline{B}_{\varepsilon'-\beta} (x'_i)\}_{i=1}^N\) are disjoint. Thus \(A' \in G(q,y,\varepsilon,k)\). 

For \(l > 1/q\), we then define
\[U(q,y,l) = \{A \in \mathcal{K}(X): A \subseteq (\overline{B}_{q-1/l}(y))^c\}.\]
As \((\overline{B}_{q-1/l}(y))^c\) is open, \(U(q,y,l)\) is open. We observe that
\[S =  \bigcap_{q \in \mathbb{Q}^+} \bigcap_{y \in D}\left(\left(\bigcap_{i=1}^{\infty} \bigcap_{k=\lfloor 1/s \rfloor+1 }^{\infty}G\left(q,y, \frac{1}{i}, k \right)\right)\cup \bigcap_{l =\lfloor 1/q \rfloor +1}^{\infty}U(q,y,l)\right),\]
so \(S\) is \(G_{\delta}\). 
\end{proof}
\begin{proposition}
\label{Prop4}
    (\cite{robinsonSharples2013})  Let \(X\) and \(Y\) be metric spaces. For compact subsets \(F \subset X\) and \(G \subset Y\) the upper box dimension of the product set \(F \times G\) satisfies
    \[\overline{\dim}_B (F) + \underline{\dim}_B(G) \leq \overline{\dim}_B (F \times G) \leq \overline{\dim}_B(F) + \overline{\dim}_B (G).\]
\end{proposition}
We first establish the packing-dimension estimate when $X$ is zero-dimensional. In this case, the approximation provides disjoint closed covers, allowing us to construct finite unions of product sets with the required local box dimension bounds.
\begin{proposition}
\label{7.4}
    Let \(X\) be a compact zero-dimensional metric space with \(\mathrm{sl-}\dim_B(X) < \infty\). A typical \(A \in \mathcal{C} (X,[0,1]^k)\) satisfies
    \begin{align*}
        \dim_P (A )\geq \mathrm{sl-}\overline{\dim}_B(A) =\mathrm{sl-}\overline{\dim}_B(X)+k.
    \end{align*}
\end{proposition}
\begin{proof}
The first relation (the inequality) follows from Lemma \ref{Lem7.3}, so it is enough to prove the second relation (the equality). We first show \(\mathrm{sl-}\overline{\dim}_B(A) \leq \mathrm{sl-}\overline{\dim}_B(X )+k\). For each \(i \in \mathbb{N}\), pick \(r_i\) and \(x_i \in X \) such that
\[\overline{\dim}_B (\overline{B}_{r_i}(x_i)) < \mathrm{sl-}\overline{\dim}_B(X) + \frac{1}{i}. \]
As \(A\) has full projection, there exists \(z_i \in A\) with \(\mathrm{proj}_1(z_i) = x_i\). Hence, we conclude by Proposition \ref{Prop4}
\[\mathrm{sl-}\overline{\dim}_B (A) \leq \overline{\dim}_B(\overline{B}_{r_i}(z_i) \cap A)\leq \overline{\dim}_B(\overline{B}_{r_i}(x_i)\times [0,1]^k)
< \mathrm{sl-}\overline{\dim}_B(X)+k + \frac{1}{i}.\]
As \(i \in \mathbb{N}\) was arbitrary, we conclude
\[\mathrm{sl-}\overline{\dim}_B (A) \leq \mathrm{sl-}\overline{\dim}_B(X )+k .\]

Now let \(G = \{A \in \mathcal{C} (X,[0,1]^k): \mathrm{sl-}\overline{\dim}_B(A) = \mathrm{sl-}\overline{\dim}_B(X)+k\}\). By Lemma \ref{Lem7.4}, \(G\) is \(G_{\delta}\). Thus it suffices to show that \(G\) is dense.

Fix \(E \in \mathcal{C} (X,[0,1]^k)\) and \(\varepsilon > 0\). By Lemma \ref{2.7} with (c), we obtain a finite closed \(\varepsilon/2\)-cover \(\{F_i\}_{i=1}^N\) of \(X\), a finite open \(\varepsilon/2\)-cover \(\{V_j\}_{j=1}^L\) of \(Y\), and an index set \(\mathcal{J} \subset [N] \times [L]\) with the following properties
\begin{enumerate}[(1)]
    \item \(\mathrm{proj}_1 (\mathcal{J}) = [N]\).
    \item If we choose nonempty closed subsets \(D(i,j) \subset F_i \times V_j\) for each \((i,j) \in \mathcal{J}\) such that for every \(i \in [N]\) there exists a \(j_i\) with \(\mathrm{proj}_1 (D(i,j_i)) = F_i\) then the set \(E' = \bigcup_{(i,j) \in \mathcal{J}} D(i,j)\) is in \(\mathcal{C} (X,Y)\) and \(\d_H (E,E') \leq \varepsilon\).
    \item The cover \(\{F_i\}_{i\in [N]}\) is disjoint.
\end{enumerate}
As all \(V_j\) are open, let \(r\) be sufficiently small such that there exists \(y_j \in V_j\) for each \(j \in [L]\) such that \(\overline{B}_r(y_j) \subseteq V_j\). We then simply let \(D(i,j) \coloneq F_i \times \overline{B}_r(y_j) \) for every \((i,j) \in \mathcal{J}\) and let \(E' \coloneq \bigcup_{(i,j) \in \mathcal{J}} D(i,j)\). By (2), \(E' \in \mathcal{C} (X,Y)\) and \(\d_H (E,E') \leq \varepsilon\). By definition of \(\mathrm{sl-}\overline{\dim}_B\) and Proposition \ref{Prop4}, \(\mathrm{sl-}\overline{\dim}_B(F_i\times \overline{B}_s(y)) = \mathrm{sl-}\overline{\dim}_B(F_i)+k \) for all \(y \in Y\) and \(s \in (0,\infty)\). Thus, by (1) and (3), we conclude that
\[\mathrm{sl-}\overline{\dim}_B(E')\geq \mathrm{sl-}\overline{\dim}_B(X)+k.\]
Hence \(E' \in G\) as desired.
\end{proof}
To extend the preceding result to general compact metric spaces, we reduce to the zero-dimensional case by finding compact subsets of $X$ whose packing dimensions approach that of $X$.
The following measure-theoretic results allow us to carry out this reduction.
\begin{definition}
    In a metric space, let \(\mathcal{P}^s\) denote the radius-based packing measure as defined e.~g. in \cite{bishop_peres_2017}. 
\end{definition}
\begin{remark}
    Proposition 2.7.1 of \cite{bishop_peres_2017} attests that this radius-based packing measure induces the standard definition of packing dimension equal to modified upper box dimension. We are using the radius-based definition because in \cite{JoycePreiss1995} the following Proposition \ref{Prop1} is proved only for that definition. In fact, it is proved in \cite{Joyce1997} that in some metric spaces Proposition \ref{Prop1} would be false with the diameter-based definition of packing measure.
\end{remark}
\begin{proposition}
\label{Prop1}(\cite{JoycePreiss1995}) Suppose \(X\) is an analytic metric space such that \(\mathcal{P}^s (X)=  \infty\) for some \(s \in (0,\infty)\). Then \(X\) has a compact subset of positive finite \(\mathcal{P}^s\) measure. 
\end{proposition}
\begin{proposition}
\label{Prop2}(\cite{Głodkowski2024}) Suppose that \(K\) is a metrizable compact space and \(\mu\) is a non-zero Radon measure on \(K\). Then there is a compact zero-dimensional subset \(Z \subseteq K\) such that \(\mu(Z) >  0\).
\end{proposition}
\begin{proposition}
\label{Prop3}(\cite{BishopPeres1996}) Let \(E\) be a subset of a separable metric space, with \(\dim_P(E) > \alpha\). Then there is a relatively closed nonempty subset \(\tilde{E}\) of \(E\), such that \(\dim_P(\tilde{E} \cap V ) > \alpha\) for any open set \(V\) which intersects \(\tilde{E}\).
\end{proposition}
\begin{theorem}
\label{7main}
Let \(X\) be a compact metric space with \(\dim_P (X) \in (0,\infty)\). A typical \(A \in \mathcal{C} (X,[0,1]^k)\) satisfies
\[\dim_P(A)= \dim_P (X )+k.\]
\end{theorem}
\begin{proof}
Let \(\dim_P(X) = s\). By Proposition \ref{Prop1}, for every \(n \in \mathbb{N}\), \(n > 2/s\), there exists compact \(K_n \subseteq X\) with positive finite \(\mathcal{P}^{s-1/n}\) measure. Then, \(\mathcal{P}^{s-1/n} \) on \(K_n\) is a finite Borel measure on a Polish space, so is therefore a Radon measure. Therefore, by Proposition \ref{Prop2}, there exists a compact zero-dimensional \(Z_n \subseteq K_n\) with \(\mathcal{P}^{s-1/n} (Z_n) > 0\). Finally, as \(\dim_P(Z_n) > s-2/n\), we conclude by Proposition \ref{Prop3} that there exists compact \(H_n \subseteq Z_n\) with \(\mathrm{sl-}\overline{\dim}_B(H_n) \geq s -2 /n\). By Lemma \ref{Lem7.3}, \(\mathrm{sl-}\overline{\dim}_B(H_n) < \infty\), so we now apply Proposition \ref{7.4} and conclude that for typical \(A \in \mathcal{C} (H_n,[0,1]^k)\)
\[\dim_P(A) \geq \mathrm{sl-}\overline{\dim}_B (H_n \times [0,1]^k) \geq s -2/n + k.\]
As \(\mathrm{sl-}\overline{\dim}_B (H_n) > 0\), \(H_n\) is perfect, so by Lemma \ref{restriction}, we conclude that for typical \(A \in \mathcal{C} (X,[0,1]^k)\)
\[\dim_P(A) \geq \dim_P(A \cap (H_n \times [0,1]^k)) \geq s-2/n + k.\]
As \(n\in \mathbb{N}\) was arbitrary, by the Baire category theorem, we conclude that for typical \(A \in \mathcal{C} (X,[0,1]^k)\), \(\dim_P(A) = s+k\).
\end{proof}
We now apply the packing-dimension estimates to collections of affine lines. Passing from a direction and translation vector to the corresponding affine line forgets the component of translation along that line, and can therefore decrease packing dimension by at most one.
\begin{definition}
    Let \(\mathrm{Graff}_1 (\mathbb{R}^n)\) be the space of affine lines, let \(\pi: \mathbb{RP}^{n-1} \times \mathbb{R}^n \to \mathrm{Graff}_1(\mathbb{R}^n)\) be given by \(\pi (l,v) = l + v\), and let \(\pi_1 : \mathrm{Graff}_1 (\mathbb{R}^n) \to \mathbb{RP}^{n-1}\) be the canonical projection sending an affine line to its unique parallel linear line.
\end{definition}
\begin{theorem}
\label{7.12}
     For all \(B \subseteq \mathbb{RP}^{n-1} \times \mathbb{R}^n\), we have that
    \[\dim_P(\pi(B)) \geq \dim_P(B) -1\]
\end{theorem}
\begin{proof}
    Let \(L \coloneq \pi(B)\). For every \(l \in L\), we will find an open neighborhood \(U_l \subseteq L\) such that there exists a locally bilipschitz homeomorphism \(\varphi_l: \pi^{-1} (U_l) \to U_l \times \mathbb{R}\) with \(\mathrm{proj}_1 \circ \varphi_l = \pi\).  For sufficiently small \(\varepsilon > 0\), we may associate a unit vector \(v_{\lambda} \in S^{n-1}\) to each \(\lambda \in \pi_1(B_{\varepsilon}(l))\)  such that \(\lambda =\mathrm{span}(v_{\lambda})\) and the map \(\lambda \mapsto v_{\lambda}\) is a bilipschitz homeomorphism. Let \(U_l \coloneq B_{\varepsilon}(l) \cap L\). To each \(l' \in U_l\), we associate a unique point \(p(l')\) on \(l'\) such that \(p(l') \perp v_{\pi_1(l')}\). Then for every \(x \in l'\), we have that \(x-p(l') = c(l',x) \cdot v_{\pi_1(l')}\) for some \(c(l',x) \in \mathbb{R}\). Thus, we may define a function \(\varphi_l:  \pi^{-1} (U_l) \to U_l \times \mathbb{R}\) given by
\[\varphi_l(\lambda , x) = (\lambda+x, c(\lambda+x, x)).\]
As \(p(\lambda+x)\) and \(v_{\lambda}\) are locally Lipschitz functions of \(\lambda\) and \(x\) on \(\pi^{-1}(U_l)\), \(\varphi_l\) is clearly locally Lipschitz.

Meanwhile, \(\varphi_l\) also admits a locally Lipschitz inverse
\[\varphi_l^{-1} (l', y) = (l'- (p(l') + y v_{\pi_1(l')}), p(l') + y v_{\pi_1(l')}),\]
so clearly \(\varphi_l\) is a locally bilipschitz homeomorphism.

Now as \(L\) is Lindelöf, \(\{U_l\}_{l \in L}\) admits a countable cover \(\{U_{l_i}\}_{i=1}^{\infty}\) of \(L = \pi(B)\). As each \(\varphi_{l_i}\) is locally bilipschitz and the packing dimension is countably stable, by the product formula for packing dimension (see \cite{Howroyd1996}), we have that
    \begin{align*}
        \dim_P(B) \leq\dim_P (\pi^{-1}(L)) =\sup_{i \in \mathbb{N}}(\dim_P(\pi^{-1}(U_{l_i}))) &= \sup_{i \in \mathbb{N}} (\dim_P (U_{l_i} \times \mathbb{R})) \\
        &\leq \sup_{i \in \mathbb{N}}(\dim_P(U_{l_i}) + 1) \\
        & = \dim_P(\pi(B)) + 1.
    \end{align*}
\end{proof}
\begin{corollary}
\label{7.13}
A typical Besicovitch set \(K \subset \mathbb{R}^n\) has the following property. If \(L\) is a set of lines that contains at least one line in each direction such that \(l \cap K\) contains a unit interval for each \(l \in L\), then
\[\dim_P(L) = 2n-2.\]
\end{corollary}
\begin{proof}
Since \(\dim_P(\mathrm{Graff}_1(\mathbb{R}^n)) = 2n-2\), it is enough to prove \(\dim_P(L) \geq 2n-2\). Let \(K_C(A) \subset \mathbb{R}^n\) be a Besicovitch set which is coded by \(A \in \mathcal{C} (\mathbb{RP}^{n-1},[0,1]^n)\) and let \(L\) be given as above. If \(\pi: \mathbb{RP}^{n-1} \times \mathbb{R}^n \to \mathrm{Graff}_1 (\mathbb{R}^n)\) is given by \(\pi(l,v) = l+v\), we may observe that typically \(L \subseteq \pi(A)\) by Corollary \ref{good1}.

Take some Cantor set \(H \subset \mathbb{RP}^{n-1}\) with \(\dim_P(H) = n-1\), which one may construct using a standard fat Cantor set construction. Now let \(B \coloneq \pi^{-1} (L) \cap A\). By definition, \(\mathrm{proj}_1 (B) = \pi_1 (\pi(B)) = \mathbb{RP}^{n-1}\). However, as \(H\) is a Cantor set, by Theorem \ref{3.3}, Lemma \ref{restriction}, and the Baire category theorem, a typical code set \(A \in \mathcal{C} (\mathbb{RP}^{n-1},[0,1]^n)\) satisfies \(\# (A \cap (\{x\} \times [0,1]^n)) = 1\) for each \(x \in H\). Thus \(A \cap (H \times [0,1]^n) \subseteq B\). By Lemma \ref{restriction} and Theorem \ref{7main}, we further conclude that typically
\begin{align*}\dim_P (B) \geq \dim_P(A \cap (H \times [0,1]^n) ) \geq \dim_P(H) + n = 2n-1. \end{align*}
Finally, we apply Theorem \ref{7.12} to conclude that
\[\dim_P(L) \geq  \dim_P(B) -1 = 2n-2.\]
\end{proof}
\begin{corollary}
\label{8.19}
    For any \(n > 1\), a typical Besicovitch set \(K \subset \mathbb{R}^n\) is not sticky. 
\end{corollary}
\begin{proof}
Immediate from Corollary \ref{7.13}.
\end{proof}
\section{Typical Code Sets are Purely Unrectifiable}
We proceed after Tam\'as K\'atay (\cite{katay2025generic}, who proved the case \(k=n=1\) for all directions other than the horizontal) and show that typical code sets \(A \in \mathcal{C}_k([0,1]^n)\) have orthogonal projections of measure zero in almost all directions, which will suffice to conclude that typical \(A \in \mathcal{C}_k([0,1]^n)\) are purely unrectifiable.
\begin{definition}
Let \(K\) be a self-similar set generated by the iterated function system \(\{\psi_i: X \to X| i \in [m]\}\). We say that \(\{\psi_i\}_{i\in [m]}\) (resp. \(K\)) satisfies the open set condition if there exists nonempty open \(V \subseteq X\) such that
\begin{itemize}
    \item \(\bigcup_{i=1}^m  \psi_i(V) \subseteq V\)
    \item The sets \(\psi_1(V), \dots, \psi_m(V)\) are pairwise disjoint.
\end{itemize}
\end{definition}
\begin{lemma}
\label{6.3}(\cite{Mattila1982StructureSelfSimilarFractals}) If \(K \subseteq \mathbb{R}^n\) is a self-similar \(k\)-set with the open set condition and \(K\) is not contained in any \(k\)-dimensional affine space then \(K\) is purely \(k\)-unrectifiable.
\end{lemma}
We use this criterion to construct compact purely unrectifiable sets with full projection onto the first coordinate and finite positive Hausdorff measure.
These sets form the basis for the subsequent density argument.
\begin{lemma}
\label{5.3}
Let \(X = \mathbb{R}^k\) and \(Y = \mathbb{R}^n\). For any axis-parallel cube \(Q \subseteq X \times Y\), there exists a compact purely \(k\)-unrectifiable \(A \subseteq Q\) such that \(\mathrm{proj}_1(A) = \mathrm{proj}_1 (Q)\) and \(\mathcal{H}^k(A) \in (0,\infty)\).
\end{lemma}
\begin{proof}
We assume without loss of generality that \(Q = [0,1]^k \times [0,1]^n\) and scale and translate as needed. Let \(S \coloneq \{0,1/3,2/3\}^k\). We construct a function \(f: S \to [0,2/3]^n\) such that \(f(0) = 0\), and if \(G(f)\) denotes the graph of \(f\),
\[\dim(\mathrm{span}(G(f))) \geq k+1.\]
First we let \(f(0) = 0\) as required. We then let \(\{e_i\}_{i=1}^k\) be the canonical basis of \(X\). We have that \(\left\{\frac{1}{3} e_i\right\}_{i=1}^k \subseteq S\). We impose
\[f\left(\frac{1}{3} e_i\right) = 0\]
for all \(i \in [k]\) and \(f(2/3, 0, \dots, 0) = (1/3, 0, \dots, 0) \neq 0\). We define \(f\) on the remainder of \(S\) arbitrarily.

We have that \(\left(\frac{1}{3}e_i,0\right) \in G(f)\) for all \(i \in [k]\), so 
\[\mathrm{span}(G(f)) \supseteq \mathrm{span} \left(\left\{\left(\frac{1}{3} e_i,0\right)\right\}_{i=1}^k\right) = X \times \{0\}.\]
We further have that \(((2/3, 0,\dots, 0), (1/3, 0, \dots, 0)) \in G(f) \cap (X \times \{0\})^c\), so \(\mathrm{span}(G(f)) \supseteq X \times \{0\}\) but \(\mathrm{span}(G(f)) \neq X \times \{0\}\). Thus \(\dim (\mathrm{span} (G(f))) > \dim (X \times \{0\}) = k\). 

We now let \(T \coloneq G(f)\). We note that \(T \subseteq [0,2/3]^{k+n}\) and \(\mathrm{proj}_1: T \to S\) is a bijection. We define mappings \(\psi_t: \mathbb{R}^{k+n} \to \mathbb{R}^{k+n}\) for each \(t \in T\) as follows
\begin{align*}
    \psi_t(A)= \frac{1}{3} A + t.
\end{align*}
We then define a function \(\Psi: \mathcal{K}(\mathbb{R}^{k+n}) \to \mathcal{K}(\mathbb{R}^{k+n})\) given by
\[\Psi (A) = \bigcup_{t \in T} \psi_t (A) \]
and let \(K\) be the unique fixed point of this iterated function system. It is easy to check that the \(\{\psi_t\}_{t \in T}\) satisfy the open set condition for open set \(V = (0,1)^{k+n}\). Further, \(K\) is a \(k\)-set with \(\mathcal{H}^k(K) \in (0,\infty)\) as each \(\psi_t\) scales by \(1/3\) and \(\# (T)= 3^k\). 

Now we must check that \(K\) is not contained in any affine \(k\)-space. However, as \(0 \in T\), \(0 \in \Psi^i([0,1]^{k+n})\) for all \(i \in \mathbb{N}\), so we conclude that
\[\Psi^i([0,1]^{k+n}) \supseteq \Psi(0) = T.\]
Thus \(T \subseteq K\) so it suffices to show that \(T\) is not contained in any affine \(k\)-space. As \(0 \in T\), any affine \(k\)-space containing \(T\) must be a linear space. But by construction, \(\dim (\mathrm{span} (T)) \geq k+1\), so \(T\) cannot be contained in any linear \(k\)-space. Thus, as \(K\) is a self-similar \(k\)-set with the open set condition, \(K\) is purely unrectifiable by Lemma \ref{6.3}.

Finally, we show that \(\mathrm{proj}_1(K) = [0,1]^k\). We write \(K\) as a nested intersection of compact sets
\[K = \bigcap_{i=1}^{\infty}\Psi^i([0,1]^{k+n})\]
We also see that \(\mathrm{proj}_1 (T) = \{0,1/3,2/3\}^k\) which suffices to show that \(\mathrm{proj}_1 (A) = [0,1]^k \implies \mathrm{proj}_1 (\Psi(A)) = [0,1]^k\). Thus, as the \(\Psi^i([0,1]^{k+n})\) are compact and nested, we conclude that
\[\mathrm{proj}_1(K) = \mathrm{proj}_1 \left(\bigcap_{i=1}^{\infty}\Psi^i([0,1]^{k+n})\right) = \bigcap_{i=1}^{\infty}\mathrm{proj}_1 (\Psi^i([0,1]^{k+n}))= [0,1]^k.\]
\end{proof}
The preceding construction allows us to approximate arbitrary code sets by finite unions of purely unrectifiable sets with finite $k$-dimensional Hausdorff measure.
By the Besicovitch-Federer projection theorem, these approximations have projections of measure zero in almost every $k$-dimensional direction.
\begin{theorem}
\label{5.5}
    A typical \(A \in \mathcal{C}_k([0,1]^n)\) has orthogonal projections of measure zero onto almost every \(L \in \mathbf{Gr}_k(\mathbb{R}^{k+n})\).
\end{theorem}
\begin{proof}
    We seek to show that
    \[G = \left\{A \in \mathcal{C}_k([0,1]^n): \text{a.e. } L \in \mathbf{Gr}_k(\mathbb{R}^{k+n}), \lambda_k(\mathrm{proj}_L (A)) = 0\right\}\]
    is a dense \(G_{\delta}\) set. We write \(G = \bigcap_{i=1}^\infty G_i\) with
    \[G_i = \left\{A \in \mathcal{C}_k ([0,1]^n): \mu \left(\left\{L \in \mathbf{Gr}_k(\mathbb{R}^{k+n}): \lambda_k (\mathrm{proj}_L(A)) \geq \frac{1}{i}\right\}\right) < \frac{1}{i}\right\}\]
    where \(\mu\) is the natural measure on the Grassmannian and we claim that each \(G_i\) is open.

    Take arbitrary \(A \in G_i\) and fix \(L \in \mathbf{Gr}_k(\mathbb{R}^{k+n})\) such that \(\lambda_k(\mathrm{proj}_L(A)) < 1/i\). We claim that there exists sufficiently large \(m \in \mathbb{N}\) such that 
    \[\lambda_k \left(\mathrm{proj}_L \left(\overline{B}^{1/m} (A)\right)\right) < \frac{1}{i}.\]
    Clearly as \(A\) is compact \(A = \bigcap_{m=1}^{\infty}\overline{B}^{1/m} (A)\). 
    As all \(\overline{B}^{1/m} (A)\) are compact, we have \(\lambda_k(\mathrm{proj}_L(\overline{B}^{1/m} (A)) )< \infty\), so
    \[\lim_{m \to \infty} \lambda_k\left(\mathrm{proj}_L\left(\overline{B}^{1/m} (A)\right)\right) = \lambda_k \left(\bigcap_{m=1}^{\infty}\mathrm{proj}_L \left( \overline{B}^{1/m} (A)\right)\right) = \lambda_k(\mathrm{proj}_L (A)) < \frac{1}{i}.\]
    This suffices to conclude the claim. 
    
    However, we now claim that there exists sufficiently large \(m \in \mathbb{N}\) such that \(\overline{B}^{1/m} (A) \cap [0,1]^{k+n}\in G_i\), which implies that \(G_i\) is open. From the prior claim, we have that
    \begin{align*}&\bigcap_{m=1}^\infty \left\{L \in \mathbf{Gr}_k (\mathbb{R}^{k+n}) : \lambda_k \left(\mathrm{proj}_L\left(\overline{B}^{1/m}(A)\right)\right) \geq \frac{1}{i}\right\} \\
    = & \left\{L \in \mathbf{Gr}_k (\mathbb{R}^{k+n}): \lambda_k (\mathrm{proj}_L(A)) \geq \frac{1}{i}\right\}.\end{align*}
    By definition of the natural measure \(\mu\), \(\mu(\mathbf{Gr}_k(\mathbb{R}^{k+n})) < \infty\), so we conclude
    \begin{align*}&\lim_{m \to \infty} \mu\left(\left\{L \in \mathbf{Gr}_k (\mathbb{R}^{k+n}) : \lambda_k \left(\mathrm{proj}_L\left(\overline{B}^{1/m}(A)\right)\right) \geq \frac{1}{i}\right\}\right) \\
    = & \mu\left(\left\{L \in \mathbf{Gr}_k (\mathbb{R}^{k+n}): \lambda_k (\mathrm{proj}_L(A)) \geq \frac{1}{i}\right\}\right) < \frac{1}{i}.\end{align*}
    Thus there exists sufficiently large \(m \in \mathbb{N}\) with \(\overline{B}^{1/m} (A) \cap [0,1]^{k+n}\in G_i\), so \(G_i\) is open.

    It therefore suffices to show that \(G\) is dense. Take \(E \in \mathcal{C}_k([0,1]^n) \) and \(\varepsilon > 0\). As \(E\) is compact, cover \(E\) with a finite collection of axis-parallel cubes \(\{Q_j\}_{j=1}^N \subseteq [0,1]^{k+n}\) with \(\diam (Q_j) < \varepsilon\) and all \(Q_j \cap E \neq \emptyset\). By Lemma \ref{5.3}, there exists a compact purely unrectifiable \(A_j \subseteq Q_j\)  with \(\mathrm{proj}_1 (A_j) = \mathrm{proj}_1 (Q_j)\). Define
    \[E' = \bigcup_{j=1}^{N} A_j.\]
    Clearly \(E'\) is compact and \(\mathrm{proj}_1 (E') = [0,1]^k\), so \(E' \in \mathcal{C}_k ([0,1]^n)\). Further, \(E'\) is a finite union of purely unrectifiable sets of measure \(\mathcal{H}^k(E') \in (0,\infty)\) so \(E'\) is purely \(k\)-unrectifiable with measure \(\mathcal{H}^k(E') \in (0,\infty)\). Thus, by the Besicovitch-Federer projection theorem, \(E' \in G\). Also, \(\d_H(E, E') \leq \max_i\diam (Q_i) < \varepsilon\) by construction, so \(G\) is dense in \(\mathcal{C}_k([0,1]^n)\). 
\end{proof}
\begin{corollary}
\label{9.5} A typical code set \(A \in \mathcal{C}_k ([0,1]^n)\) is purely \(k\)-unrectifiable.
\end{corollary}
\begin{proof}
Immediate from Proposition \ref{5.2}, Theorem \ref{5.5}, and the Besicovitch-Federer projection theorem.
\end{proof}
\begin{conjecture}
    A typical \(A \in \mathcal{C}_k ([0,1]^n)\) has orthogonal projections of measure zero onto every \(L \in \mathbf{Gr}_k(\mathbb{R}^{k+n}) \setminus \{\mathrm{span}(X)\}\).
\end{conjecture}
\section{Typical Code Sets Have Small Injective Projections}
We now study projections that are simultaneously small in Hausdorff dimension and injective on typical code sets. The key idea is to construct sets that retain their full first-coordinate projections while collapsing to single points under a fixed projection. By separating these images along a line, we obtain approximations with small projection fibers, letting us establish typical injectivity.
\begin{lemma}
\label{8.2}
Let \(X = [0,1]^k\) and \(Y = [0,1]^n\). If \(L\) is a linear subspace of \(\mathrm{span}(X) \times \mathrm{span}(Y)\) with \(\mathrm{proj}_1 (L^{\perp}) = \mathrm{span}(X)\) and \(l \subseteq L\) is a line through the origin, then there exists a constant \(\rho > 0\) with the following property. If \(Q_1 \subseteq \mathrm{span}(X)\) and \(Q_2 \subseteq \mathrm{span}(Y)\) are closed axis parallel cubes with \(\diam Q_1 < \rho \diam Q_2\), then there exists an infinite sequence of compact subsets \(\{D_j\}_{j=1}^{\infty}\subseteq Q_1 \times Q_2\) such that 
\begin{enumerate}[(a)]
    \item for each \(j\), \(\mathrm{proj}_1 (D_j) = Q_1\),
    \item for each \(j\), \(\mathrm{proj}_L(D_j)\) is a singleton,
    \item for each \(j \neq j'\), the singletons \(\mathrm{proj}_l(D_j)\) and \(\mathrm{proj}_l(D_{j'})\) are distinct.
\end{enumerate}
\end{lemma}
\begin{proof}
    Let \(Q_1'\) be an axis-parallel closed cube of \(\mathrm{span}(X)\) with diameter 1. Since \(\mathrm{proj}_1 (L^{\perp}) = \mathrm{span} (X)\), there exists a compact subset \(D\) of \(L^{\perp}\) such that \(\mathrm{proj}_1(D) = Q'_1\). Let \(Q_2'\) be an axis parallel closed cube of \(\mathrm{span}(Y)\) such that its interior contains \(\mathrm{proj}_2(D)\). Then \(D\) is a subset of \(Q_1' \times \mathrm{int} (Q_2')\) and \(\mathrm{proj}_L(D)\) is a singleton.

    Let \(\delta\) be the positive distance between \(\mathrm{proj}_2 (D)\) and the boundary of \(Q'_2\) and let \(\{v_j\}_{j=1}^{\infty} \subseteq \mathrm{span}(Y)\) be vectors of length less than \(\delta\) with distinct projections onto \(l\). Such a choice exists as otherwise \(l \subseteq \mathrm{span}(X)\), contradicting \(\mathrm{proj}_1 (L^{\perp}) = \mathrm{span}(X)\). 

    Now \(D'_j \coloneq D + v_j\) is a good choice of \(D_j\) for \(Q_1 = Q_1'\) and \(Q_2 = Q_2'\). Thus, if we choose \(\rho = 1/\diam  Q'_2\), properly scaled and translated copies of the sets \(D'_j\) as \(D_j\) have all the required properties. 
\end{proof}
The preceding construction provides the approximations needed to show that small Hausdorff dimension and injectivity of the projection hold simultaneously for typical code sets.
\begin{theorem}
\label{8.3}
Let \(X = [0,1]^k\) and \(Y = [0,1]^n\). If \(L\) is a linear subspace of \(\mathrm{span} (X) \times \mathrm{span}(Y)\) with \(L \not\subseteq \mathrm{span}(X)\), then a typical \(A \in \mathcal{C} (X,Y)\) satisfies \begin{enumerate}[(i)] \item \(\dim_H (\mathrm{proj}_L(A)) = \dim_H (L \cap (\mathrm{span}(X) \times \{0\}))\). \item \(\mathrm{proj}_L: A \to L\) is injective. \end{enumerate}
\end{theorem}
\begin{proof}
    Note that for linear subspaces \(\dim_H (L) = \dim L\). We fix an arbitrary linear subspace \(L \subseteq \mathrm{span} (X) \times \mathrm{span} (Y)\) with \(L \not\subseteq \mathrm{span} (X)\) and let \(j \coloneq \dim_H(L \cap (\mathrm{span}(X) \times \{0\}))\). 

For any \(A \in \mathcal{C} (X,Y)\), we have that
\begin{align*}\mathrm{proj}_{L \cap (\mathrm{span}(X) \times \{0\})}(\mathrm{proj}_{L} (A)) &= \mathrm{proj}_{L \cap (\mathrm{span}(X) \times \{0\})}(\mathrm{proj}_{\mathrm{span}(X) \times \{0\}}(A)) \\
&= \mathrm{proj}_{L \cap (\mathrm{span}(X) \times \{0\})}(X \times \{0\}).\end{align*}
As projection maps are open, \(\mathrm{proj}_{L \cap (\mathrm{span}(X) \times \{0\})}(X \times \{0\})\) has nonempty interior in \(L \cap (\mathrm{span}(X) \times \{0\})\) as \(X \times \{0\}\) has nonempty interior in \(\mathrm{span}(X) \times \{0\}\). Thus
\[\dim_H(\mathrm{proj}_{L \cap (\mathrm{span}(X) \times \{0\})}(X \times \{0\})) = \dim_H(L \cap (\mathrm{span}(X) \times \{0\})) = j. \]
Therefore, as projections are Lipschitz,
\[\dim_H(\mathrm{proj}_L(A)) \geq \dim_H (\mathrm{proj}_{L \cap (\mathrm{span}(X) \times \{0\})}(\mathrm{proj}_L(A))) =  j.\]
To demonstrate (i), it thus suffices to show that a typical \(A \in \mathcal{C} (X,Y)\) has \(\dim_H(\mathrm{proj}_L(A)) \leq j\). 

Let \(G = \{A \in \mathcal{C} (X,Y): \dim_H (\mathrm{proj}_L(A)) \leq j\}\). By Lemma \ref{Lem2.1}, \(G\) is a \(G_{\delta}\) set. Thus to show (i), it suffices to show that \(G\) is dense. 

We now construct a subspace \(\hat{L}\) of \(L\) such that \(\dim \hat{L} = \dim L - j\) and \(\mathrm{proj}_1(\hat{L}^{\perp}) = \mathrm{span}(X)\). We begin by letting \(\Lambda \) be the (possibly trivial) orthogonal complement in \(\mathrm{span}(X)\) of \(\mathrm{proj}_1 (L^{\perp})\). Then
\begin{align*}\dim \Lambda = k -  \dim (\mathrm{proj}_1(L^{\perp})) &= \dim(\mathrm{span}(X)) -\dim(\mathrm{proj}_{L^{\perp}}(\mathrm{span}(X) \times \{0\})) \\
& = \dim (L \cap (\mathrm{span}(X) \times \{0\})) = j\end{align*}
where in the penultimate equality we used the rank-nullity theorem for \(\mathrm{proj}_{L^{\perp}} |_{\mathrm{span}(X) \times \{0\}}\). We then let \(\hat{L} \coloneq (L^{\perp} + \Lambda \times \{0\})^{\perp}\), that is the orthogonal complement of  \(\Lambda \times \{0\}\) in \(L\). Now
\[\mathrm{proj}_1 (\hat{L}^{\perp}) = \mathrm{proj}_1 (L^{\perp} + \Lambda \times \{0\}) = \mathrm{proj}_1 (L^{\perp}) + \Lambda = \mathrm{span}(X)\]
as desired. Further, as \(\dim \Lambda = j\), \(\dim \hat{L} = \dim L - j\).  

Next, we define a \(G_{\delta}\) set corresponding to condition (ii). Pick any line \(l \subseteq \hat{L} \subseteq L\). We have that \(\mathrm{proj}_l : A \to l\) factors through \(\mathrm{proj}_L: A \to L\), so if \(\mathrm{proj}_l: A \to l\) is injective, that suffices to show that \(\mathrm{proj}_L: A \to L\) is injective. Let \(\{H_z\}_{z \in \mathbb{R}}\) be the 1-parameter family of hyperplanes that form the fibers of \(\mathrm{proj}_l\). We have that
\begin{align*}K&\coloneq \{A \in \mathcal{C}(X,Y)|\mathrm{proj}_l: A \to l \text{ is an injection}\} \\
&=\left\{A \in \mathcal{C} (X,Y): \forall z \in \mathbb{R}, \# (A \cap H_z) \leq 1\right\}.\end{align*}
Thus it suffices to show that \(K\) is a dense \(G_{\delta}\) set. We decompose \(K\) into a countable intersection of sets \(K = \bigcap_{m=1}^{\infty} K_m\), with
\[K_m = \left\{A \in \mathcal{C} (X,Y): \forall z \in \mathbb{R}, \forall p_1,p_2 \in A \cap H_z,  \d_{X \times Y} (p_1,p_2) < \frac{1}{m}\right\}.\]
If we let \(\mathcal{D}= \{H_z \cap(X \times Y) \text{ nonempty}:z \in \mathbb{R}\}\), \(K_m = G_m(\mathcal{D},2)\) in the notation of Lemma \ref{2.3}, so all \(K_m\) are open. Thus by the Baire category theorem it suffices to show that each \(K_m\) is dense.

We now take arbitrary \(E \in \mathcal{C}(X,Y)\) and \(\varepsilon \in (0,1/m)\). We will construct \(E' \in G \cap K_m\) with \(\d_H (E,E') \leq \varepsilon\). We apply Lemma \ref{8.2} with respect to \(\hat{L}\) and \(l\) to obtain \(\rho > 0\) with the property that if \(Q_1 \subseteq \mathrm{span}(X)\) and \(Q_2 \subseteq \mathrm{span}(Y)\) are closed axis-parallel cubes with \(\diam  Q_1 < \rho \diam Q_2\), then there exists an infinite sequence of compact subsets \(\{D_j\}_{j=1}^{\infty}\) such that
    \begin{enumerate}[(a)]
    \item for each \(j\), \(\mathrm{proj}_1 (D_j) = Q_1\),
    \item for each \(j\), \(\mathrm{proj}_{\hat{L}}(D_j)\) is a singleton,
    \item for each \(j \neq j'\), the singletons \(\mathrm{proj}_l(D_j)\) and \(\mathrm{proj}_l(D_{j'})\) are distinct.
\end{enumerate}
We now cover \(E\) with axis parallel cuboids \(\{\tilde{Q}_i\}_{i=1}^N \subseteq [0,1]^{k+n}\), \(\tilde{Q}_i \coloneq Q_{i,1} \times Q_{i,2}\) satisfying \(\diam Q_{i,1} < \rho \diam Q_{i,2}\), \(\diam \tilde{Q}_i \leq \varepsilon\), and all \(\tilde{Q}_i \cap E \neq \emptyset\). Let \(\{D_{i,j}\}_{j=1}^{\infty} \subseteq \tilde{Q}_i\) denote the sequence of compact subsets given by the lemma. As each collection is infinite and satisfies condition (c), inductively, we may choose a \(\tilde{D}_i \in \{D_{i,j}\}_{j=1}^{\infty}\) for each \(i \in [N]\) such that all \(\mathrm{proj}_l(\tilde{D}_i)\) are distinct. Define
\[E' = \bigcup_{i=1}^N \tilde{D}_i.\]
Clearly \(E'\) is compact, and as by condition (a), \(\mathrm{proj}_1(\tilde{D}_i) = Q_{i,1}\) for each \(i \in [N]\),
\[\mathrm{proj}_1(E') = \bigcup_{i=1}^N Q_{i,1} \supseteq \mathrm{proj}_1 (E) = X.\]
Thus \(E' \in \mathcal{C} (X,Y)\). Further, as every \(x \in E\) and \(x' \in E'\) is contained in some \(\tilde{Q}_i\), we have that \(\d_H (E,E') \leq \max_i\diam (\tilde{Q}_i) \leq \varepsilon\).

Thus, to show that \(G\) and \(K_m\) are dense, it suffices to show that \(E' \in G \cap K_m\). First we show that \(E' \in G\). By condition (b), each \(\mathrm{proj}_{\hat{L}} (\tilde{D}_i)\) is a singleton, so \(\# (\mathrm{proj}_{\hat{L}} (E')) < \infty\). As \(\Lambda\) is the orthogonal complement of \(\hat{L}\) in \(L\), \(\mathrm{proj}_{L \to \hat{L}}: L\to \hat{L}\) has kernel \(\Lambda \), where \(\mathrm{proj}_{L \to \hat{L}}\) just denotes the restriction of \(\mathrm{proj}_{\hat{L}}\) to \(L\). Further, \(\mathrm{proj}_{L \to \hat{L}} (\mathrm{proj}_L(E')) = \mathrm{proj}_{\hat{L}}(E')\), so
\[\mathrm{proj}_L (E') \subseteq \mathrm{proj}^{-1}_{L \to \hat{L}} (\mathrm{proj}_{\hat{L}}(E')) = \bigcup_{x \in \mathrm{proj}_{\hat{L}}(E')} (\Lambda + x)\]
as each fiber of the projection is a shifted copy of \(\Lambda\). However, \(\dim_H(\Lambda) = \dim \Lambda = j\) and \(\mathrm{proj}_{\hat{L}}(E')\) is finite, so we conclude that \(\dim_H(\mathrm{proj}_L(E')) \leq j\) so \(E' \in G\). 

Finally, we show that \(E' \in K_m\). By construction of the \(\tilde{D}_i\), each \(\mathrm{proj}_l (\tilde{D}_i)\) is a distinct singleton. Thus \(p_1,p_2 \in E'\cap H_z\) for some \(z \in \mathbb{R}\) if and only if \(\mathrm{proj}_l(p_1) = \mathrm{proj}_l (p_2)\) if and only if \(p_1,p_2 \in \tilde{D}_i\) for some fixed \(i \in [N]\). However \(\tilde{D}_i \subseteq \tilde{Q}_i\) so \(\diam (\tilde{D}_i) \leq \varepsilon\). Thus whenever \(p_1,p_2 \in E'\cap H_z\) for some \(z\), \(\d_{X \times Y}(p_1,p_2) \leq \varepsilon < 1/m\), so \(E' \in K_m\).
\end{proof}
\begin{corollary}
\label{10.3}
    Let \(X = [0,1]^k\) and \(Y = [0,1]^n\) and let \(l \subseteq \mathrm{span}(Y)\) be a line through the origin. If we let \(L \coloneq \mathrm{span}(X) \times l \), then for typical \(A \in \mathcal{C}(X,Y)\), \(\mathrm{proj}_L(A)\) is the graph of a continuous surjection \(f: S\to X\) for some \(S \subseteq l\).
\end{corollary}
\begin{proof}
    Fix \(l \subseteq \mathrm{span}(Y)\). Let \(\{H_z\}_{z \in l}\) denote the 1-parameter family of hyperplanes perpendicular to \(l\) and take \(z \in l\) as the parameter such that \(l \cap H_{z} = z\). By Theorem \ref{8.3}, a typical \(A \in \mathcal{C} (X,Y)\) satisfies \(\# (A \cap H_z) \leq 1\) for all \(z \in l\) as the \(H_z\) are the fibers of \(\mathrm{proj}_l\). Projecting onto \(L\), we have that 
    \[\#(\mathrm{proj}_L(A) \cap \mathrm{proj}_L(H_z)) \leq 1.\]
    As projections onto nested subspaces commute, the hyperplanes \(\mathrm{proj}_L(H_z)\) remain the fibers of \(\mathrm{proj}_l : L  \to l\). Thus let \(S \coloneq \{z \in l:\#(\mathrm{proj}_L(A) \cap \mathrm{proj}_L(H_z)) = 1\}\). We define \(f: S \to X\) such that \(f(z)\) is the \(x\)-coordinate of the singleton element of \(\mathrm{proj}_L(A) \cap \mathrm{proj}_L(H_z)\). Then \(\mathrm{proj}_L(A)\) is a graph of \(f\). However, 
    \[\mathrm{proj}_1 (\mathrm{proj}_L(A)) = \mathrm{proj}_1(A) = X\]
    so \(f\) surjects onto \(X\). As \(X\) is compact and the graph of \(f\) is closed, \(f\) is continuous.
\end{proof}
We conclude by applying the projection theorem to obtain restrictions on the intersections of partial Besicovitch sets with horizontal hyperplanes.
\begin{corollary}
\label{8.5}
    Let \(H \subseteq \mathbb{R}^{k+1}\) be any fixed horizontal hyperplane and let \(l \subseteq H\) be an arbitrary line. Then let \(\{h_z\}_{z \in \mathbb{R}} \subseteq H\) be the 1-parameter family of affine \(k-1\)-dimensional hyperplanes of \(H\) perpendicular to \(l\). For a typical partial Besicovitch set \(K \subseteq \mathbb{R}^{k+1}\)
    \begin{enumerate}[label=(\arabic*)]
        \item \(\dim_H(K \cap H) = 0\),
        \item for any \(z \in \mathbb{R}\), \(K\) contains at most one line intersecting \(h_z\).
    \end{enumerate}
\end{corollary}
\begin{proof}
    Let a partial Besicovitch set \(K \subseteq \mathbb{R}^{k+1}\) be coded by \(A \in \mathcal{C}_k([0,1]^k)\) and let the horizontal hyperplane \(H\) be given by \(H \coloneq \mathbb{R}^k \times \{z_0\}\). We first show (1). A point \((u,z_0)\) is contained in the line coded by \((m,b)\) if and only if \(u = mz_0 + b\), where \(u,m,b \in \mathbb{R}^k\) and \(z_0 \in \mathbb{R}\). As Corollary \ref{good1} tells us there are no uncoded lines in \(K\), we see that \((u,z_0) \in H \cap K\) if and only if \((m,b) \in \{(x, u-xz_0):x \in  \mathbb{R}^k\}\) for some \((m,b) \in A\). Let \(L \subseteq \mathbb{R}^{2k}\) be the linear subspace 
    \[L \coloneq (\{(x,-xz_0): x \in \mathbb{R}^k\})^{\perp}  = \begin{cases}
    \mathrm{span}(Y), & z_0 = 0 \\
    \{(x,x/z_0): x \in \mathbb{R}^k\}, & \text{otherwise}.
    \end{cases}\]
    We see that \(\mathrm{proj}_{L} (m,b) = \mathrm{proj}_{L}(0,u)\), so in general \(\mathrm{proj}_L\{(0,u) \in \mathbb{R}^{2k}: (u,z_0) \in H \cap K\}= \mathrm{proj}_{L}(A)\). Thus, as \(\mathrm{proj}_L\) is a bilipschitz homeomorphism on \(\{0\} \times \mathbb{R}^k\)
    \[\dim_H (H \cap K) \leq \dim_H(\mathrm{proj}_{L}(A)).\]
    Clearly \(L \cap (\mathrm{span}(X) \times \{0\}) = 0\), so by Theorem \ref{8.3} (i), \(\dim_H(\mathrm{proj}_{L}(A)) = 0\) for typical \(A \in \mathcal{C}_k([0,1]^k)\). Thus \(\dim_H(H \cap K) = 0\) for typical partial Kakeya sets \(K\). 
    
    Now we show (2). Let \(v_0 \in \mathbb{R}^k\) be the unit direction vector of \(l\). The one parameter family of \(k-1\)-dimensional hyperplanes \(\{h_z\}_{z \in \mathbb{R}}\) is then of the form \(h_c \coloneq \{(u,z_0): \langle u, v_0 \rangle = c\}\) where \(c \in \mathbb{R}\). The line coded by \((m,b)\) intersects \(h_c\) if and only if \(\langle mz_0+b, v_0 \rangle = c\), which can be rewritten as \(\langle (m,b) , (z_0 v_0,v_0) \rangle = c\). The set of these code points forms an affine hyperplane \(H_c\) in \(\mathbb{R}^{2k}\) perpendicular to the fixed vector \((z_0 v_0, v_0)\) which is not contained in \(\mathrm{span}(X)\) as \(v_0\) is a nonzero vector. We now apply Theorem \ref{8.3} (ii) with respect to \(l = \mathrm{span}(z_0v_0,v_0)\) and conclude that for typical \(A \in \mathcal{C}_k([0,1]^k)\), \(\mathrm{proj}_l :A \to l\) is injective. Thus, as the hyperplanes \(H_c\) are the fibers of this projection \(\mathrm{proj}_l\), we conclude that \(\#(H_c \cap A)\leq 1\) for all \(c\). Therefore, for a typical partial Kakeya set \(K\), \(K\) contains at most one line intersecting \(h_c\) for each \(c\). 
\end{proof}
\medskip

\printbibliography

\end{document}